\documentclass[12pt]{article}
\usepackage{amsthm,amsmath,amssymb,amscd,verbatim}
\usepackage{graphics,pifont}
\usepackage{amssymb}
\usepackage{graphicx}
\date{}
\newcommand{\al}{\alpha}
\newcommand{\be}{\beta}
\newcommand{\ga}{\gamma}

\newcommand{\sq}{$\square$}
\newcommand{\va}{\varepsilon}
\begin{document}
\title{On the Chaos in Continuous Weakly Mixing Maps}                       
\author{Bau-Sen Du\footnote {Bau-Sen Du is a retired research fellow at the Institute of Mathematics, Academia Sinica, Taiwan} \\ 
dubs@gate.sinica.edu.tw \cr}
\maketitle

\begin{abstract}
Let $\mathcal X$ be an infinite locally compact separable metric space with metric $\rho$ and let $f : \mathcal X \longrightarrow \mathcal X$ be a continuous weakly mixing map.  Let $\beta = \sup \big\{ \rho(x, y): \{x, y \} \subset \mathcal X \big\}$.  In this note, we use a method which generalizes the classical construction of the Cantor ternary set in $[0, 1]$ to show (Theorem 4) that, for any countably infinite set $\{x_1, x_2, \cdots\}$ of points in $\mathcal X$ with compact orbit closures $\overline{O_f(x_i)}$'s, there exist an infinite set $\mathcal M$ of positive integers and countably infinitely many pairwise disjoint Cantor sets ${\mathcal S}^{(1)}, {\mathcal S}^{(2)}, \cdots$ of totally transitive points of $f$ in $\mathcal X$ such that (1) for any integers $\ell \ge 1$ and $n \ge 1$, $\ell!$ divides all sufficiently large integers in $\mathcal M$ and for any distinct points $a_1, a_2, \cdots, a_n$ in ${\mathbb S} = \bigcup_{j=1}^\infty \, {\mathcal S}^{(j)}$, the set $\{ F_n^m\big((a_1, a_2, \cdots, a_n)\big): m \in \mathcal M \}$ is dense in $\mathcal X \times \mathcal X \times \cdots \times \mathcal X$ ($n$ terms), where $F_n\big((a_1, a_2, \cdots, a_n)\big) = \big(f(a_1), f(a_2), \cdots, f(a_n)\big)$; (2) ${\mathbb S}$ is a dense $\beta$-scrambled set of $f^n$ for all $n \ge 1$; (3) for any $x$ in $\{x_1, x_2, \cdots\}$ and any $c$ in $\widehat {\mathbb S} = \bigcup_{i=0}^\infty \, f^i({\mathbb S})$, the set $\{ x, c \}$ is a (${\beta}/2$)-scrambled set of $f$.  Furthermore, if $f$ has a fixed point and $\delta = \inf_{n \ge 1} \big\{ \sup\{ \rho(f^n(x), x): x \in \mathcal X \} \big\} \ge 0$, then the above Cantor sets ${\mathcal S}^{(1)}, {\mathcal S}^{(2)}, \cdots$ can be chosen to satisfy the additional property that $\widehat {\mathbb S} = \bigcup_{i=0}^\infty f^i({\mathbb S})$ is a dense {\it invariant} $\delta$-scrambled set of $f^n$ for all $n \ge 1$.  For continuous mixing maps on $\mathcal X$, we have a similar but stronger result (Theorem 5) which is used to show the existence of dense {\it invariant} $\xi$-scrambled sets, where $\xi > 0$, for continuous transitive maps on an interval (Corollary 6) although such transitive maps may not be mixing.  A notion of chaos is also introduced. \vspace*{.1in}


\end{abstract}

The problem on the existence of dense uncountable scrambled ({\it invariant} or not) sets for continuous (topologically) weakly mixing maps on some nontrivial dynamical systems have long been addressed by many people, e.g., Babilonov\'a {\bf\cite{bab}}, Balibrea, Guirao and Oprocha {\bf\cite{bali}}, Boro\'nski, Kupka and Oprocha {\bf\cite{boro}}, Fory$\acute {\rm s}$, Huang, Li and Oprocha {\bf\cite{fo}}, Iwanik {\bf\cite{iw}}, Li and Ye {\bf\cite{liye}}. Ruette {\bf\cite{ru}}, Xiong and Yang {\bf\cite{xiong}} and Yuan and L$\ddot{{\rm u}}$ {\bf\cite{yuan}} and so on.  In this note, we re-address this problem by a different method which generalizes the classical construction of the Cantor ternary set in the unit interval $[0, 1]$ to any continuous (topologically) weakly mixing map on some metric space in such way that at each stage we take ever more steps to construct sequences of finite unions of nested compact sets whose total sums of diameters tend to zero.  Since our method works for any infinite locally compact separable metric space, in the sequel, we let $(\mathcal X, \rho)$ denote any infinite locally compact separable metric space with metric $\rho$.  

Recall that, since $\mathcal X$ is separable, it has a countable dense subset and, since $\mathcal X$ is locally compact, for each point $x$ in $\mathcal X$, there is an open neighborhood $U$ of $x$ such that the closure $\overline U$ is compact.  For infinite separable locally compact metric spaces, we have the following well-known result.

\noindent
{\bf Lemma 1.}
{\it Let $\mathcal X$ be an infinite locally compact separable metric space.  Then, for any point $x$ in $\mathcal X$ and any open neighborhood $U$ of $x$, there is an open neighborhood $V$ of $x$ such that the closure ${\overline V}$ is a compact subset of $U$.  Furthermore, $\mathcal X$ has a countable basis of open sets with compact closures.}

Let $f : \mathcal X \longrightarrow \mathcal X$ be a continuous map and let $\mathbb N$ denote the set of all positive integers.  Let $s \ge 0$ be an integer.  For any point $x_0$ and any nonempty open sets $U$ and $V$ in $\mathcal X$, we define \\
[8pt]\centerline{$N(x_0, V) = \{ n \in \mathbb N: f^n(x_0) \in V \}$ and $N\big(f^s(U), V\big) = \{ n \in \mathbb N: f^{n+s}(U) \cap V \ne \emptyset \}$.} \\
[8pt]\indent We say that $f$ is (topologically) transitive if, for any nonempty open sets $W_1$ and $W_2$ in $\mathcal X$, there is a positive integer $n$ such that $W_1 \cap f^{-n}(W_2) \ne \emptyset$, or, equivalently, $f^n(W_1) \cap W_2 \ne \emptyset$.  So, $n \in N(W_1, W_2) \ne \emptyset$.  In this case, $N(W_1, W_2)$ is actually an infinite subset of $\mathbb N$.  This fact will be used in the proof of Lemma 3(2) below.  We say that $f$ is (topologically) totally transitive if $f^k$ is (topologically) transitive for each $k \ge 1$. 

We say that $f$ is (topologically) weakly mixing if $f \times f: \mathcal X \times \mathcal X \longrightarrow \mathcal X \times \mathcal X$ defined by $(f \times f)\big((x, y)\big) = \big(f(x), f(y)\big)$, is (topologically) transitive, i.e ., for any nonempty open sets $U_1$, $U_2$, $V_1$ and $V_2$ in $\mathcal X$, there exists a positive integer $n$ such that $(f \times f)^n\big((U_1 \times U_2)\big) \, \cap \, (V_1 \times V_2) \ne \emptyset$, i.e., $f^n(U_1) \cap V_1 \ne \emptyset$ and $f^n(U_2) \cap V_2 \ne \emptyset$, in particular, $f$ is (topologically) transitive.   

For continuous weakly mixing maps, we have the following well-known equivalent statements in which Part (2) is the celebrated Furstenberg Intersection Lemma {\bf\cite{fur}}.  

\noindent
{\bf Lemma 2.}
{\it Let $(\mathcal X, \rho)$ be an infinite locally compact separable metric space with metric $\rho$ and let $f : \mathcal X \rightarrow \mathcal X$ be a continuous weakly mixing map.  Then the following statements are equivalent: 
\begin{itemize}
\item[{\rm (1)}]
$f$ is continuous weakly mixing, i.e., $f \times f: \mathcal X \times \mathcal X \longrightarrow \mathcal X \times \mathcal X$ is transitive;  

\item[{\rm (2)}]
for any integer $n \ge 2$ and any nonempty open sets $U_j, V_j$, $1 \le j \le n$, in $\mathcal X$, there exist nonempty open sets $\widehat U_n$ and $\widehat V_n$ in $\mathcal X$ such that 
$$
\emptyset \ne N(\widehat U_n, \widehat V_n) \subset N(U_1, V_1) \cap N(U_2, V_2) \cap N(U_3, V_3) \cap \cdots \cap N(U_n, V_n).
$$
Consequently, the map $f \times f \times \cdots \times f (n$ terms) : $\mathcal X \times \mathcal X \times \cdots \times \mathcal X (n$ terms) $\longrightarrow \mathcal X \times \mathcal X \times \cdots \times \mathcal X (n$ terms) is transitive;

\item[{\rm (3)}] 
for any nonempty open sets $U$ and $V$ in $\mathcal X$, the set $N(U, V)$ contains arbitrarily long string of consecutive positive integers;

\item[{\rm (4)}] 
for any nonempty open sets $U$ and $V$ in $\mathcal X$, the set $N(U, V)$ contains two consecutive positive integers;

\item[{\rm (5)}] 
$f$ is transitive and, for any nonempty open set $U$ in $\mathcal X$, the set $N(U, U)$ contains two consecutive positive integers;

\item[{\rm (6)}] 
$f$ is transitive and, for any nonempty open set $U$ in $\mathcal X$, the set $N(U, U)$ contains arbitrarily long string of consecutive positive integers.             
\end{itemize}}

\noindent
{\it Proof.}
$(1) \longrightarrow (2)$: We follow the proof of Theorem 1.51 in {\bf\cite{linear}} on p.21.  Let $n \ge 2$ be an integer and let $U_i, V_i$, $1 \le i \le n$ be nonempty open sets in $\mathcal X$.  Since $f$ is weakly mixing, $N(U_1, U_2) \cap N(V_1, V_2) \ne \emptyset$.  Let $k \in N(U_1, U_2) \cap N(V_1, V_2)$.  Then the sets $\widehat U_2 = U_1 \cap f^{-k}(U_2)$ and $\widehat V_2 = V_1 \cap f^{-k}(V_2)$ are nonempty. Since $f$ is weakly mixing, it is transitive too.  So, $N(\widehat U_2, \widehat V_2) \ne \emptyset$.  Let $m \in N(\widehat U_2, \widehat V_2)$ and let $x_0 \in \widehat U_2 = U_1 \cap f^{-k}(U_2) (\subset U_1)$ be a point such that $f^m(x_0) \in \widehat V_2 = V_1 \cap f^{-k}(V_2) (\subset V_1)$.  Then, $x_0 \in U_1$, $f^m(x_0) \in V_1$ and $f^m\big(f^k(x_0)\big) = f^k\big(f^m(x_0)\big) \in V_2$.  Since $x_0 \in \widehat U_2 \subset f^{-k}(U_2)$, $f^k(x_0) \in U_2$.  Therefore, we obtain that $N(\widehat U_2, \widehat V_2) \subset N(U_1, V_1) \cap N(U_2, V_2)$.  Similar arguments imply that there exist nonempty open sets $\widehat U_3$ and $\widehat V_3$ such that $N(\widehat U_3, \widehat V_3) \subset N(\widehat U_2, \widehat V_2) \cap N(U_3, V_3) \subset \big(N(U_1, V_1) \cap N(U_2, V_2)\big) \cap N(U_3, V_3)$.  Now it follows easily from induction that there exist nonempty open sets $\widehat U_n$ and $\widehat V_n$ in $\mathcal X$ such that $N(\widehat U_n, \widehat V_n) \subset N(U_1, V_1) \cap N(U_2, V_2) \cap N(U_3, V_3) \cap \cdots N(U_n, V_n)$.

$(2) \longrightarrow (3)$: For any integer $n \ge 2$ and any nonempty open sets $U$ and $V$ in $\mathcal X$, it follows from Part (2) that there exist nonempty open sets $U_n'$ and $V_n'$ in $\mathcal X$ such that 
$$N(U_n', V_n') \subset N\big(U, f^{-1}(V)\big) \cap N\big(U, f^{-2}(V)\big) \cap N\big(U, f^{-3}(V)\big) \cap \cdots \cap N\big(U, f^{-n}(V)\big).$$  
It is easy to see that, for any integer $m$ in $N(U_n', V_n')$, $\{ m+1, m+2, m+3, \cdots, m+n \} \subset N(U, V)$.  Therefore, $N(U, V)$ contains $n$ consecutive positive integers for each $n \ge 2$.

$(5) \longrightarrow (6)$: The following proof is due to Huang and Ye {\bf\cite{hy}}.  Let $U$ be any nonempty open set.  Assume that $N(U, U)$ contains $k \ge 2$ consecutive positive integers $n, n+1, \cdots, n+k-1$.  So, \\
[8pt]\centerline{$U \cap f^{-n}(U) \ne \emptyset, \, U \cap f^{-n-1}(U) \ne \emptyset, \, \cdots, \, U \cap f^{-n-k+1}(U) \ne \emptyset$.} \\
[8pt]\indent Since $f$ is transitive, there exist positive integers $j_1, j_2, \cdots, j_{k-1}$ such that the sets \\
[8pt]\hspace*{1.1321in} $U \cap f^{-n}(U) \ne \emptyset$, \\
[5pt]\hspace*{1.12in} $\big(U \cap f^{-n}(U)\big) \cap f^{-j_1}(U \cap f^{-n-1}(U)) \ne \emptyset$, \\
[5pt]\centerline{$\big((U \cap f^{-n}(U)) \cap f^{-j_1}(U \cap f^{-n-1}(U))\big) \cap f^{-j_2}(U \cap f^{-n-2}(U)) \ne \emptyset$,} \\
[1pt]\centerline{\vdots \qquad\qquad\qquad\qquad\qquad \vdots} \\
[4pt]and $\big(U \cap f^{-n}(U)\big) \cap f^{-j_1}\big(U \cap f^{-n-1}(U)\big) \cap f^{-j_2}\big(U \cap f^{-n-i}(U) \cap \cdots \cap f^{-j_{k-1}}\big(U \cap f^{-n-k+1}(U)\big) \ne \emptyset$. \\
[8pt]\indent Let $W = \big(U \cap f^{-n}(U)\big) \cap f^{-j_1}\big(U \cap f^{-n-1}(U)\big) \cap f^{-j_2}\big(U \cap f^{-n-i}(U) \cap \cdots \cap f^{-j_{k-1}}\big(U \cap f^{-n-k+1}(U)\big)$.  Then $W \ne \emptyset$.  By assumption, there exist two consecutive positive integers $m$ and $m+1$ such that \\
$W \cap f^{-m}(W) \ne \emptyset$ and $W \cap f^{-m-1}(W) \ne \emptyset$.  So, we have $\emptyset \ne W \cap f^{-m}(W) \subset \big(U \cap f^{-n}(U)\big) \cap f^{-m}\big(U \cap f^{-n}(U)\big) \subset U \cap f^{-m-n}(U)$.  On the other hand, for each $1 \le i \le k-1$, we consider \\
$\emptyset \ne W \subset \big(U \cap f^{-n}(U)\big) \cap f^{-j_i}\big(U \cap f^{-n-i}(U)\big)$ and obtain that \\
[5pt]\hspace*{.45in}$\emptyset \ne W \cap f^{-m}(W)$ $\subset \big(f^{-j_i}(U \cap f^{-n-i}(U))\big) \,\, \cap \,\, f^{-m}\big(f^{-j_i}(U \cap f^{-n-i}(U))\big)$ \\
[5pt]\hspace*{1.767in}$\subset f^{-j_i}(U) \cap f^{-j_i}\big(f^{-m-n-i}(U)\big) = f^{-j_i}\big(U \cap f^{-m-n-i}(U)\big)$ and \\
[5pt]\hspace*{.45in}$\emptyset \ne W \cap f^{-m-1}(W) \subset \big(f^{-j_i}(U \cap f^{-n-i}(U))\big) \,\, \cap \,\, f^{-m-1}\big(f^{-j_i}(U \cap f^{-n-i}(U))\big)$ \\
[5pt]\hspace*{1.915in}$\subset f^{-j_i}(U) \, \cap \, f^{-j_i}(f^{-m-1}(f^{-n-i}(U))) = f^{-j_i}\big(U \cap f^{-m-n-1-i}(U)\big)$. \\
[8pt]\indent In particular, we obtain that \\
[8pt]\centerline{$U \cap f^{-m-n}(U) \ne \emptyset$, $U \cap f^{-m-n-1}(U) \ne \emptyset$, $U \cap f^{-m-n-2}(U) \ne \emptyset$, $\cdots$, and $U \cap f^{-m-n-k}(U) \ne \emptyset$,} \\
[8pt]i.e., $N(U, U)$ contains the $k+1$ consecutive positive integers $m+n$, $m+n+1$, $m+n+2$, $\cdots$, and $m+n+k$.  

$(6) \longrightarrow (1)$: Let $U_1, U_2, V_1, V_2$ be any nonempty open sets in $\mathcal X$.  Since $f$ is transitive, there exist positive integers $i, j, k$ successively such that $i < j < k$ and $k > j-i > 0$ and \\
[8pt]\centerline{$U_1 \cap f^{-i}(U_2) \ne \emptyset, \big(U_1 \cap f^{-i}(U_2)\big) \cap f^{-j}(V_2) \ne \emptyset$, and $\big((U_1 \cap f^{-i}(U_2)) \cap f^{-j}(V_2)\big) \cap f^{-k}(V_1) \ne \emptyset$.}\\
[8pt]\indent Let $W = \big(\big(U_1 \cap f^{-i}(U_2)\big) \cap f^{-j}(V_2)\big) \cap f^{-k}(V_1)$.  Then $W \ne \emptyset$.  By assumption, $N(W, W)$ contains arbitrarily long string of consecutive positive integers, say, $n, n+1, n+2, \cdots, n+k+i-j, \cdots, n+k-1, n+k$.  In particular, $n \in N(W, W)$.  Let $w$ be a point in $W$ $(\subset U_1)$ such that $f^n(w) \in W$ $(\subset f^{-k}(V_1))$.  So, $f^{n+k}(w) \in V_1$, i.e., $n+k \in N(U_1, V_1)$.  Furthermore, $N(W,W)$ contains $n+k+i-j$.  Let $w'$ be a point in $W$ $(\subset U_1 \cap f^{-i}(U_2))$ such that $f^{n+k+i-j}(w') \in W$ $(\subset f^{-j}(V_2))$.  Then $w' \in U_1$, $f^i(w') \in U_2$ and $f^{n+k}(f^i(w')) = f^{n+k+i}(w') \in f^j(W) \subset V_2$.  So, $n+k \in N(U_2, V_2)$.  This, combined with the above, implies that $n+k \in N(U_1, V_1) \cap N(U_2, V_2)$.  Therefore, $f$ is weakly mixing.
\hfill\sq
 
We shall base our method of constructing the appropriate scrambled sets on the following result in which Part (1) is an easy consequence of Lemma 2(3) while Part (2) is essentially due to Akin and Kolyada (see the proof of Theorem 3.8 in {\bf\cite{akin}}).

\noindent
{\bf Lemma 3.}
{\it Let $\mathcal X$ be an infinite locally compact separable metric space and let $f : \mathcal X \rightarrow \mathcal X$ be a continuous weakly mixing map.  Then the following hold:
\begin{itemize}
\item[{\rm (1)}]
For any integer $n \ge 2$, let  $K_1, \, K_2, \, \cdots, \, K_n$ be compact sets in $\mathcal X$ with nonempty interiors and let $V_1, \, V_2, \, \cdots, \, V_n$ be nonempty open sets in $\mathcal X$.  Then there exist arbitrarily long string of {\rm consecutive} positive integers $k$ and corresponding compact sets $K_{k,j} \subset K_j$, $1 \le j \le n$, with nonempty interiors such that, for each such $k$, $f^{k}(K_{k,j}) \subset V_j$ for all $1 \le j \le n$.  

\item[{\rm (2)}]
Suppose $f$ has a point $x_0$ in $\mathcal X$ whose orbit $O_f(x_0)$ has compact closure $\overline{O_f(x_0)}$.  Let $s \ge 0$ be an integer and let $\mathcal G(x_0)$ be an open cover of $\overline{O_f(x_0)}$.  Then there exists an (open) member, say $G_s(x_0)$, of $\mathcal G(x_0)$ such that, for any integer $n \ge 2$, any compact sets $K_j$, $1 \le j \le n$, with nonempty interiors, and any nonempty open sets $V_j$, $1 \le j \le n$, in $\mathcal X$, by considering the sets 
$$K_1, \, K_2, \, \cdots, \, K_n: \, V_1, \, V_2, \, \cdots, \, V_n \, ; \, G_s(x_0) \,\,\, \text{\big(note this} \,\,\, G_s(x_0) \,\,\, \text{at the end}\big),$$
there exist infinitly many positive integers $k_i = k_{i,n,s}$, $i \ge 1$ and corresponding compact sets $K_{k_i,j} \subset K_j$, $1 \le j \le n$, with nonempty interiors such that, 
$$\text{for each} \,\,\, i \ge 1, \, f^{k_i}\big(f^s(K_{k_i,j})\big) \subset V_j \,\,\, \text{for all} \,\,\, 1 \le j \le n \,\,\, \text{and} \,\,\, f^{k_i}(x_0) \in G_s(x_0).$$
\end{itemize}}

\noindent
{\bf Remarks.}
(1) In Lemma 3(1), since there exist arbitrarily long string of {\rm consecutive} positive integers $k$ and corresponding compact sets $K_{k,j} \subset K_j$, $1 \le j \le n$, with nonempty interiors such that, for each such $k$, $f^{k}(K_{k,j}) \subset V_j$ for all $1 \le j \le n$, we can choose, for any integers $\ell \ge 1$ and $s \ge 0$, in any one such long string of integers of length $> \ell! + s$, an integer $k'+s$ such that $\ell!$ divides $k'$ and $f^{k'}\big(f^s(K_{k'+s,j})\big) = f^{k'+s}(K_{k'+s,j}) \subset V_j$ for all $1 \le j \le n$.

(2) In Lemma 3(2), since we only require the sets $K_{k_i,j}$'s to satisfy $f^{k_i}\big(f^s(K_{k_i,j})\big) \subset V_j$ for all $1 \le j \le n$, we can choose them as small as we wish.

\begin{proof}
Part (1) follows from Parts (2) $\&$ (3) of Lemma 2 trivially.  We now prove Part (2).  Since $\mathcal G(x_0)$ is an open cover of the compact set $\overline{O_f(x_0)}$, it has a finite subcover, say, $\{ G_1, G_2, \cdots, G_m \}$.  So, we have \\
[8pt]\centerline{$N(x_0, G_1) \cup N(x_0, G_2) \cup N(x_0, G_3) \cup \cdots \cup N(x_0, G_m) = \mathbb N$, the set of all positive integers.} \\
[8pt]\indent Recall that $s \ge 0$ is a given integer.  Suppose, for each $1 \le i \le m$, there exists nonempty open sets $W_{i}$ and $\widehat W_{i}$ in $\mathcal X$ such that $N\big(f^s(W_{i}), \widehat W_{i}\big) \cap N(x_0, G_{i})$ is a finite set of positive integers.  Then, the finite union $\bigcup_{{i}=1}^m \, \big(N(f^s(W_{i}), \widehat W_{i}) \, \cap \, N(x_0, G_{i})\big)$ is finite.  Since
$\bigcup_{{i}=1}^m \, \big(N(f^s(W_{i}), \widehat W_{i}) \, \cap \, N(x_0, G_{i})\big) \supset \bigcup_{{i}=1}^m \, \big((\bigcap_{\ell=1}^m N(f^s(W_{\ell}), \widehat W_{\ell})) \, \bigcap \, N(x_0, G_{i})\big) = \big(\bigcap_{\ell=1}^m N(f^s(W_{\ell}), \widehat W_{\ell})\big) \cap \big(\bigcup_{{i}=1}^m N(x_0, G_{i})\big) =$ \\
$\bigcap_{\ell=1}^m N\big(f^s(W_{\ell}), \widehat W_{\ell}\big)$,
this implies that $\bigcap_{\ell=1}^m N\big(f^s(W_{\ell}), \widehat W_{\ell}\big)$ is a finite set of positive integers.  Since \\
[8pt]\centerline{$k \in \bigcap_{\ell=1}^m N\big(W_{\ell}, \widehat W_{\ell}\big)$ and $k > s$ \, if and only if \, $k-s \in \bigcap_{\ell =1}^m N \big(f^s(W_{\ell}), \widehat W_{\ell}\big)$,} \\  
[8pt]we obtain that $\bigcap_{\ell=1}^m N\big(W_{\ell}, \widehat W_{\ell}\big)$ is a finite set of positive integers. 

However, it follows from Lemma 2(2) that there exist nonempty open sets $\widetilde W$ and $\widetilde W'$ such that $N(\widetilde W, \widetilde W') \subset \bigcap_{\ell=1}^m N\big(W_\ell, \widehat W_\ell\big)$ and so, $N\big(\widetilde W, \widetilde W'\big)$ is a finite set of positive integers which contradicts the fact that $f$ is weakly mixing and so transitive.  Therefore, there exists one (open) member, say $G_s(x_0)$, of $\mathcal G(x_0)$ such that, \\
[8pt]\hspace*{.5in}for any nonempty open sets $W$ and $W'$ in $\mathcal X$, \\
[5pt]\hspace*{1.5in}$N\big(f^s(W), W'\big) \cap N\big(x_0, G_s(x_0)\big)$ is an infinite set of positive integers. \\
[8pt]\indent For any $n \ge 2$, any compact sets $K_1, K_2, \cdots, K_n$ with nonempty interiors and any nonempty open sets $V_1, V_2, \cdots, V_n$, by considering the nonempty {\it interiors} $K_j^\circ$ of $K_j$ and nonempty open sets $V_j$, it follows from Lemma 2(2) that there exist nonempty open sets $U$ and $V$ such that $N(U, V) \subset \bigcap_{j=1}^n N\big(K_j^\circ, V_j\big)$.  Since, if $k \in N\big(f^s(U), V\big)$, then $k+s \in N(U, V) \subset \bigcap_{j=1}^n N\big(K_j^\circ, V_j\big)$ which implies that $k \in \bigcap_{j=1}^n N\big(f^s(K_j^\circ), V_j\big)$, this shows that $N\big(f^s(U), V\big) \subset \bigcap_{j=1}^n N\big(f^s(K_j^\circ), V_j\big)$. 

Since $N\big(f^s(U), V\big) \cap N\big(x_0, G_s(x_0)\big)$ is an infinite set of positive integers, so is the larger set \\
$\big(\bigcap_{j=1}^n N\big(f^s(K_j^\circ), V_j\big)\big) \cap N\big(x_0, G_s(x_0)\big)$.  Let $<k_i \,(= k_{i,n,s})>_{i \ge 1}$ be any infinite sequence of {\it such} positive integers.  Then, for any $i \ge 1$ and all $1 \le j \le n$, $k_i \in N(f^s(K_j^\circ), V_j) \cap N(x_0, G_s(x_0))$.  For each $1 \le j \le n$, let $c_{i,j}$ be a point in $K_j^\circ$ such that $f^{k_i+s}(c_{i,j})$ $\in V_j$.  By Lemma 1, there is an open neighborhood ${\mathbb W}_j$ of $f^{k_i+s}(c_{i,j})$ such that the closure $\overline{\mathbb W_j}$ is compact and $\overline{\mathbb W_j} \subset V_j$.  For each $1 \le j \le n$, let $E_{k_i,j}$ be an open neighborhood of $c_{i,j}$ in $K_j^\circ$ such that $f^{k_i+s}(E_{k_i,j}) \subset \mathbb W_j$.  Let $K_{k_i, j} = \overline{E_{k_i,j}}$.  Then $K_{k_i, j}$ is a compact subset of $K_j$ such that $K_{k_i,j}^\circ \supset E_{k_i,j} \ne \emptyset$.  Therefore, we have \\
[8pt]\centerline {$f^{k_i}\big(f^s(K_{k_i,j})\big) = f^{k_i+s}(K_{k_i,j}) = f^{k_i+s}(\overline{E_{k_i,j}}) = \overline{f^{k_i+s}(E_{k_i,j})} \subset \overline{\mathbb W_j} \subset V_j$.}\\
[8pt]This shows that there exist infinitely many positive integers $k_i \,(= k_{i,n,s})$, $i \ge 1$, and corresponding compact sets $K_{k_i,j} \subset K_j$ in $\mathcal X$ with nonempty interiors $K_{k_i,j}^\circ \supset E_{k_i,j}$, $1 \le j \le n$, such that, for each $i \ge 1$, $f^{k_i}\big(f^s(K_{k_i,j})\big) \subset V_j$ for all $1 \le j \le n$ and $f^{k_i}(x_0) \in G_s(x_0)$.\qedhere
\end{proof}

Let $\mathcal M$ be an infinite set of positive integers.  Let $C$ be a nonempty set and $v$ a point in $\mathcal X$.  Following Mai {\bf\cite{mai}}, we say that the set $C$ is synchronously proximal ($\mathcal M$-synchronously proximal respectively) to the point $v$ if there exists a strictly increasing sequence $< m_i >_{i \ge 1}$ of positive integers (in $\mathcal M$ respectively) such that the diameters of the sets $f^{m_i}(C) \cup \{ v \}$ tend to zero as $m_i$ tends to $\infty$.  On the other hand, we say that $C$ is {\it dynamically} synchronously proximal to the point $v$ if there exists a strictly increasing sequence $< n_i >_{i \ge 1}$ of positive integers such that the diameters of the sets $f^{n_i}\big(C \cup \{ v \}\big)$ tend to zero as $n_i$ tends to $\infty$.  It is possible that $C$ is ({\it dynamically} respectively) synchronously proximal to many different points.  We say that $C$ is a Cantor set if it is compact, perfect and totally disconnected.

Let $S$ be a subset of $\mathcal X$ with at least two distinct points and let $\eta$ be a positive number.  We say that $S$ is a $0$-scrambled set of $f$ if, for any two distinct points $x$ and $y$ in $S$, we have 
$$\liminf_{n \to \infty} \rho\big(f^n(x), f^n(y)\big) = 0 \,\,\, \text{and} \,\, \limsup_{n \to \infty} \rho\big(f^n(x), f^n(y)\big) > 0 \,\, \big(\text{note that it is} \, > 0, \, \text{not} \, \ge 0\big).
$$
We say that $S$ is an $\eta$-scrambled set of $f$ if, for any two distinct points $x$ and $y$ in $S$, we have 
$$\liminf_{n \to \infty} \rho\big(f^n(x), f^n(y)\big) = 0 \quad \text{and} \quad \limsup_{n \to \infty} \rho\big(f^n(x), f^n(y)\big) \ge \eta.
$$  
We say that $S$ is an $\infty$-scrambled set of $f$ if, for any two distinct points $x$ and $y$ in $S$, we have 
$$\liminf_{n \to \infty} \rho\big(f^n(x), f^n(y)\big) = 0 \quad \text{and} \quad \limsup_{n \to \infty} \rho\big(f^n(x), f^n(y)\big) = \infty.
$$
\indent Let $\beta = \eta$ or $\beta = \infty$.  Inspired by Xiong and Yang {\bf\cite{xiong}}, we call $S$ a $\beta$-scrambled set (with respect to $\mathcal M$) of $f$ if, for any two distinct points $x$ and $y$ in $S$, we have
$$
\liminf_{\substack{n \to \infty \\ n \in \mathcal M}} \rho\big(f^n(x), f^n(y)\big) = 0 \,\,\, \text{and} \,\,\, \limsup_{\substack{n \to \infty \\ n \in \mathcal M}} \rho\big(f^n(x), f^n(y)\big) \ge \beta.
$$
\indent In the following result, Part (1) is a generalization of a result of Ruette {\bf\cite{ru}}, Part (2) is well-known and is an easy consequence of Part (1) while Parts (3), (4) $\&$ (5) are new.  For the proof, we shall base our arguments on a method which generalizes the classical construction of the Cantor ternary set in the unit interval $[0, 1]$.  In the classical construction, for each $\ell \ge 1$, at the $\ell^{th}$-stage, exactly one step is taken while, in our construction, we take $1+\ell+2\ell(\ell+1)+(\ell+1)^{\ell \cdot 2^{\ell+1}}=(2\ell^2 + 3\ell+1) + (\ell+1)^{\ell \cdot 2^{\ell+1}}$ steps among which \\
[6pt]\hspace*{.2in}the first step will be used to split each relevant compact set with nonempty interior into two \\
\hspace*{.5in}smaller compact subsets with nonempty interiors; \\
\hspace*{.2in}the $\ell$ steps from the $2^{nd}$ through the $\big(1+\ell\big)^{th}$ will be used to verify Parts (3) $\&$ (5), \\
\hspace*{.2in}the $2\ell(\ell+1)$ steps from the $\big((1+\ell)+1\big)^{st}$ through the $\big((1+\ell)+2\ell(\ell+1)\big)^{st}$ will be used to \\
\hspace*{.5in}confirm Part (4), and \\
\hspace*{.2in}the $(\ell+1)^{\ell \cdot 2^{\ell+1}}$ steps from the $\big((1+\ell+2\ell(\ell+1))+1\big)^{st}$ through the $\big((1+\ell+2\ell(\ell+1))+ (\ell+1)^{\ell \cdot 2^{\ell+1}}\big)^{st}$ \\
\hspace*{.5in}will be used to establish Parts (1) $\&$ (2).

\noindent
{\bf Theorem 4.}
{\it Let $(\mathcal X, \rho)$ be an infinite locally compact separable metric space with metric $\rho$ and let $f : \mathcal X \rightarrow \mathcal X$ be a continuous weakly mixing map (so, $\mathcal X$ has no isolated points).  Let $\beta$ be defined by 
$$
\beta =  \begin{cases}
               \infty, & \text{if $\mathcal X$ is unbounded}, \cr
               diam(\mathcal X) = \sup\big\{ \rho(x, y): \{ x, y \} \subset \mathcal X \big\}, & \text{if $\mathcal X$ is bounded} \cr
       \end{cases}
$$
and let $\{ v_1, v_2, v_3, \cdots \}$ be a countably infinite set of points in $\mathcal X$ and let (if any) $\{ x_1, \, x_2, \, x_3, \, \cdots \}$ be a countably infinite set of points in $\mathcal X$ whose orbits under $f$ have compact closures.

Then there exist an infinite set $\mathcal M$ of positive integers and countably infinitely many pairwise disjoint Cantor sets ${\mathcal S}^{(1)}, \, {\mathcal S}^{(2)}, \, {\mathcal S}^{(3)}, \, \cdots$ of totally transitive points of $f$ in $\mathcal X$ such that
\begin{itemize}
\item[{\rm (1)}]
for any integer $\ell \ge 1$, $\ell!$ divides all sufficiently large integers in $\mathcal M$ and, for any integer $n \ge 1$ and any distinct points $a_1, a_2, a_3, \cdots, a_n$ in $\mathbb S = \bigcup_{j=1}^\infty {\mathcal S}^{(j)}$, the set $\{F_n^{m}\big((a_1, a_2, a_3, \cdots, a_n)\big): m \in \mathcal M \}$ is dense in $\mathcal X_n = \mathcal X \times \mathcal X \times \mathcal X \times \cdots \times \mathcal X$ ($n$ terms), where $F_n: \mathcal X_n \longrightarrow \mathcal X_n$ is the map defined by $F_n\big((x_1, x_2, x_3, \cdots, x_n)\big)=\big(f(x_1), f(x_2), f(x_3), \cdots, f(x_n)\big)$.  In particular, the point $(a_1, a_2, a_3, \cdots, a_n)$ is a totally transitive point of $F_n$ for each $n \ge 1$;

\item[{\rm (2)}]
the set $\mathbb S = \bigcup_{j=1}^\infty {\mathcal S}^{(j)}$ is a dense $\beta$-scrambled set of $f^n$ in $\mathcal X$ for all $n \ge 1$; 

\item[{\rm (3)}]
for any integers $r \ge 1$ and $s \ge 0$, the set $\bigcup_{j=1}^r {\mathcal S}^{(j)}$ is $\mathcal M$-synchronously proximal to $v_m$ for all $m \ge 1$:

\item[{\rm (4)}]
if $\{x_1, x_2, x_3, \cdots\} \ne \emptyset$, then, for any integers $r \ge 1$ and $s \ge 0$, the set $f^s\big(\bigcup_{j=1}^r {\mathcal S}^{(j)}\big)$ is dynamically synchronously proximal to $x_m$ for all $m \ge 1$ and, for any point $x$ in the set $\{x_1, x_2, x_3, \cdots\}$ and any point $c$ in the set $\widehat {\mathbb S} = \bigcup_{i=0}^\infty \, f^i(\mathbb S)$, we have $\liminf_{n \to \infty} \rho\big(f^n(x), f^n(c)\big) = 0$ and $\limsup_{n \to \infty} \rho\big(f^n(x), f^n(c)\big) \ge \beta/2$, i.e., the set $\{ x, c \}$ is a ($\beta/2$)-scrambled set of $f$ in $\mathcal X$;

\item[{\rm (5)}]
if $f$ has a fixed point and $\delta = \inf_{n \ge 1} \big\{ \sup\{ \rho(f^n(x), x): x \in \mathcal X \} \big\} \ge 0$, then the above Cantor sets ${\mathcal S}^{(1)}, {\mathcal S}^{(2)}, {\mathcal S}^{(3)}, \cdots$ of totally transitive points of $f$ can be chosen (by choosing appropriate $v_1, v_2, v_3, \cdots$) so that Parts (1), (2), (3) $\&$ (4) hold and the set $\widehat {\mathbb S} = \bigcup_{i=0}^\infty \, f^i(\mathbb S)$ is a dense {\it invariant} $\delta$-scrambled set of $f^n$ in $\mathcal X$ for all $n \ge 1$.
\end{itemize}}

\noindent
{\it Proof.}
For simplicity, in the following proof, we assume that $f$ has points $x_i$'s with compact orbit closures $\overline{O_f(x_i)}$'s (such points exist when $\mathcal X$ is compact), otherwise, skip all arguments pertaining to the points $x_i$'s.  

Firstly, we introduce some notations.  

Since $\mathcal X$ is a separable metric space, it has a countable open base, say $U_1$, $U_2$, $U_3$, $\cdots$.  

Let $a_1$, $a_2$, $a_3$, $\cdots$ and $b_1$, $b_2$, $b_3$, $\cdots$ be two sequences of points in $\mathcal X$ such that
$$
\lim_{n \to \infty} \rho(a_n, b_n) = \beta = \begin{cases}
         \infty, & \text{if $\mathcal X$ is unbounded}, \\
         diam(\mathcal X) = \sup \big\{ \rho(x, y) : \{ x, y \} \subset \mathcal X \big\}, & \text{if $\mathcal X$ is bounded}. \\
         \end{cases}
$$ 
\indent For any point $x_0$ in $\mathcal X$ and any positive integer $n$, let $V_n(x_0)$ denote a nonempty open set in $\mathcal X$ of diameter $< 1/n$ which contains the point $x_0$.

For any integers $s \ge 0$, $m \ge 1$ and $k \ge 1$, let $\mathcal G_{s,k}(x_m)$ be a cover of $\overline{O_f(x_m)}$ of open sets with diameters $< 1/k$.  It follows from Lemma 3(2) that there exists a member $G_{s,k}(x_m)$ of $\mathcal G_{s,k}(x_m)$ with diameter $< 1/k$ such that, for any integer $n \ge 2$, any compact sets $K_1, K_2, \cdots, K_n$ with nonempty interiors and any nonempty open sets $V_1, V_2, \cdots, V_n$ in $\mathcal X$, the set \\
[10pt]\hspace*{.138in}$\big(N(f^s(K_1), V_1) \cap N(f^s(K_2), V_2) \cap \cdots \cap N(f^s(K_n), V_n)\big) \bigcap N\big(x_m, G_{s,k}(x_m)\big)$ \\ 
[5pt]\hspace*{4.285in}is an infinite set of positive integers. \\
[10pt]\indent For any integers $s \ge 0$, $m \ge 1$ and $k \ge 1$, let $\widehat G_{s,k}(x_m)$ be a nonempty open set in $\mathcal X$ with diameter $< 1/k$ such that the distance $dist\big(G_{s,k}(x_m), \widehat G_{s,k}(x_m)\big) = \inf\big\{ \rho(x, y): x \in G_{s,k}(x_m), y \in \widehat G_{s,k}(x_m) \big\}$ between $G_{s,k}(x_m)$ and $\widehat G_{s,k}(x_m)$ satisfies 
$$
\lim_{k \to \infty} dist\big(G_{s,k}(x_m), \widehat G_{s,k}(x_m)\big) \ge \frac {\beta}2 = \,\,\, \begin{cases}
                                                      \infty, & \text{if $\mathcal X$ is unbounded}, \\
                                                      \frac 12 diam(\mathcal X), & \text{if $\mathcal X$ is bounded}.\\
                                                      \end{cases}
$$
\indent We now start the construction which is somewhat mechanical and maybe a little boring.  However, once you get the hang of it, many steps can be omitted.

At the first stage, 

let $n_1$ be any positive integer and let $K_0^{(1,1,0)}$ and $K_1^{(1,1,0)}$ be any two disjoint compact sets in $U_1$ with nonempty interiors \big(the number 1 in the superscript $(1,*,*)$ indicates that these sets lie in $U_1$ and the number 1 in the superscript $(*,1,*)$ means that these sets are obtained at the first stage\big).

By applying Lemma 3(1) to the following $2 \cdot 2^2$ sets with nonempty interiors \\
[8pt]\centerline{$K_0^{(1,1,0)}, K_0^{(1,1,0)}, K_1^{(1,1,0)}, K_1^{(1,1,0)}: V_{n_1}(a_1), V_{n_1}(b_1), V_{n_1}(a_1), V_{n_1}(b_1)$ \big(note the subscripts of $K_{...}^{(1,1,0)}$\big),} \\
[8pt]we obtain an integer $k_{1,1} \, (> n_1)$ with $1!$ dividing $k_{1,1}$ (the number 1 in the subscript (1,*) of $k_{1,*}$ indicates that this number is obtained at the first stage) and $2^2$ pairwise disjoint compact sets $K_{00}^{(1,1,1)}$ and $K_{01}^{(1,1,1)}$ in $K_{0}^{(1,1,0)}$, \, $K_{10}^{(1,1,1)}$ and $K_{11}^{(1,1,1)}$ in $K_{1}^{(1,1,0)}$, with nonempty interiors such that \\
[10pt]\hspace*{.135in}each diameter is so smaller than $1/{2^3}$ that \\ 
[3pt]\hspace*{1.6in}the open set $U_2 \setminus \big(K_{00}^{(1,1,1)} \cup K_{01}^{(1,1,1)} \cup K_{10}^{(1,1,1)} \cup K_{11}^{(1,1,1)}\big)$ is nonempty and, \\
[10pt]\noindent for all $\al_0 = 0, 1$, $\al_1 = 0, 1$, \\
[8pt]\centerline{$f^{k_{1,1}}\big(K_{\al_0\al_1}^{(1,1,1)}\big) \subset W_{n_1}(\al_1)$, where $W_{n_1}(0) = V_{n_1}(a_1)$ and $W_{n_1}(1) = V_{n_1}(b_1)$, and $1!$ divides $k_{1,1}$.} \\
[8pt]\indent By applying Lemma 3(1) to the following $2 \cdot 2^2$ sets with nonempty interiors \\
[8pt]\centerline{$K_{00}^{(1,1,1)}, K_{01}^{(1,1,1)}, K_{10}^{(1,1,1)}, K_{11}^{(1,1,1)}: \,\, V_{n_1}(v_1), V_{n_1}(v_1), V_{n_1}(v_1), V_{n_1}(v_1)$,} \\
[8pt]we obtain an integer $k_{1,2} \, (> k_{1,1})$ with $1!$ dividing $k_{1,2}$ and $2^2$ pairwise disjoint compact sets $K_{\al_0\al_1}^{(1,1,2)} \subset K_{\al_0\al_1}^{(1,1,1)}$, $\al_0 = 0, 1$, $\al_1 = 0, 1$, with nonempty interiors such that 
$$
f^{k_{1,2}}\big(K_{00}^{(1,1,2)} \cup K_{01}^{(1,1,2)} \cup K_{10}^{(1,1,2)} \cup K_{11}^{(1,1,2)}\big) \subset V_{n_1}(v_1) \,\,\, \text{and} \,\,\, 1! \,\, \text{divides} \,\,\, k_{1,2}.
$$
\indent By applying Lemma 3(2) to the following $2 \cdot 2^2+1$ sets with nonempty interiors and $s = 0$ \\
[8pt]\centerline{$K_{00}^{(1,1,2)}, K_{01}^{(1,1,2)}, K_{10}^{(1,1,2)}, K_{11}^{(1,1,2)}: \,\, G_{0,n_1}(x_1), G_{0,n_1}(x_1), G_{0,n_1}(x_1), G_{0,n_1}(x_1); \,\, G_{0,n_1}(x_1)$,} \\
[8pt]we obtain an integer $k_{1,3} \, (> k_{1,2})$ and $2^2$ pairwise disjoint compact sets $K_{\al_0\al_1}^{(1,1,3)} \subset K_{\al_0\al_1}^{(1,1,2)}$, $\al_0 = 0, 1$, $\al_1 = 0, 1$, with nonempty interiors such that 
$$
f^{k_{1,3}}\big(K_{00}^{(1,1,3)} \cup K_{01}^{(1,1,3)} \cup K_{10}^{(1,1,3)} \cup K_{11}^{(1,1,3)}\big) \subset G_{0,n_1}(x_1) \,\,\, \text{and} \,\,\, f^{k_{1,3}}(x_1) \in G_{0,n_1}(x_1).
$$
\indent By applying Lemma 3(2) to the following $2 \cdot 2^2+1$ sets with nonempty interiors and $s = 1$ \\
[8pt]\centerline{$K_{00}^{(1,1,3)}, K_{01}^{(1,1,3)}, K_{10}^{(1,1,3)}, K_{11}^{(1,1,3)}: \,\, G_{1,n_1}(x_1), G_{1,n_1}(x_1), G_{1,n_1}(x_1), G_{1,n_1}(x_1); \,\, G_{1,n_1}(x_1)$,} \\
[8pt]we obtain an integer $k_{1,4}+1 \, (> k_{1,3}+1)$ and $2^2$ pairwise disjoint compact sets $K_{\al_0\al_1}^{(1,1,4)} \subset K_{\al_0\al_1}^{(1,1,3)}$, $\al_0 = 0, 1$, $\al_1 = 0, 1$, with nonempty interiors such that $f^{k_{1,4}+1}\big(K_{00}^{(1,1,4)} \cup K_{01}^{(1,1,4)} \cup K_{10}^{(1,1,4)} \cup K_{11}^{(1,1,4)}\big) =$
$$
f^{k_{1,4}}\big(f(K_{00}^{(1,1,4)}) \cup f(K_{01}^{(1,1,4)}) \cup f(K_{10}^{(1,1,4)}) \cup f(K_{11}^{(1,1,4)})\big) \subset G_{1,n_1}(x_1) \,\,\, \text{and} \,\,\, f^{k_{1,4}}(x_1) \in G_{1,n_1}(x_1).
$$
\indent By applying Lemma 3(2) to the following $2 \cdot 2^2+1$ sets with nonempty interiors and $s = 0$ \\
[8pt]\centerline{$K_{00}^{(1,1,4)}, K_{01}^{(1,1,4)}, K_{10}^{(1,1,4)}, K_{11}^{(1,1,4)}: \,\, \widehat G_{0,n_1}(x_1), \widehat G_{0,n_1}(x_1), \widehat G_{0,n_1}(x_1), \widehat G_{0,n_1}(x_1); \,\, G_{0,n_1}(x_1)$,} \\
[8pt]we obtain an integer $k_{1,5} \, (> k_{1,4})$ and $2^2$ pairwise disjoint compact sets $K_{\al_0\al_1}^{(1,1,5)} \subset K_{\al_0\al_1}^{(1,1,4)}$, $\al_0 = 0, 1$, $\al_1 = 0, 1$, with nonempty interiors such that 
$$
f^{k_{1,5}}\big(K_{00}^{(1,1,5)} \cup K_{01}^{(1,1,5)} \cup K_{10}^{(1,1,5)} \cup K_{11}^{(1,1,5)}\big) \subset \widehat G_{0,n_1}(x_1) \,\,\, \text{and} \,\,\, f^{k_{1,5}}(x_1) \in G_{0,n_1}(x_1).
$$
\indent By applying Lemma 3(2) to the following $2 \cdot 2^2+1$ sets with nonempty interiors and $s = 1$ \\
[8pt]\centerline{$K_{00}^{(1,1,5)}, K_{01}^{(1,1,5)}, K_{10}^{(1,1,5)}, K_{11}^{(1,1,5)}: \,\, \widehat G_{1,n_1}(x_1), \widehat G_{1,n_1}(x_1), \widehat G_{1,n_1}(x_1), \widehat G_{1,n_1}(x_1); \,\, G_{1,n_1}(x_1)$,} \\
[8pt]we obtain an integer $k_{1,6}+1 \, (> k_{1,5}+1)$ and $2^2$ pairwise disjoint compact sets $K_{\al_0\al_1}^{(1,1,6)} \subset K_{\al_0\al_1}^{(1,1,5)}$, $\al_0 = 0, 1$, $\al_1 = 0, 1$, with nonempty interiors such that $f^{k_{1,6}+1}\big(K_{00}^{(1,1,6)} \cup K_{01}^{(1,1,6)} \cup K_{10}^{(1,1,6)} \cup K_{11}^{(1,1,6)}\big) =$
$$
f^{k_{1,6}}\big(f(K_{00}^{(1,1,6)}) \cup f(K_{01}^{(1,1,6)}) \cup f(K_{10}^{(1,1,6)}) \cup f(K_{11}^{(1,1,6)})\big) \subset \widehat G_{1,n_1}(x_1) \,\,\, \text{and} \,\,\, f^{k_{1,6}}(x_1) \in G_{1,n_1}(x_1).
$$
\indent Now, there are $2^4 = 16$ finite sequences of 4 elements from the set $\{ U_1, U_2 \}$.  We label them as $W_{1,1}^{(1)}, W_{1,2}^{(1)}, W_{1,3}^{(1)}, W_{1,4}^{(1)}; W_{2,1}^{(1)}, W_{2,2}^{(1)}, W_{2,3}^{(1)}, W_{2,4}^{(1)}; W_{3,1}^{(1)}, W_{3,2}^{(1)}, W_{3,3}^{(1)}, W_{3,4}^{(1)}; \cdots$, $W_{16,1}^{(1)}$, $W_{16,2}^{(1)}$, $W_{16,3}^{(1)}$, $W_{16,4}^{(1)}$. 

\indent By applying Lemma 3(1) to the following $2 \cdot 2^2$ sets with nonempty interiors \\
[8pt]\centerline{$K_{00}^{(1,1,6)}, K_{01}^{(1,1,6)}, K_{10}^{(1,1,6)}, K_{11}^{(1,1,6)}: \,\, W_{1,1}^{(1)}, W_{1,2}^{(1)}, W_{1,3}^{(1)}, W_{1,4}^{(1)}$,} \\ 
[8pt]we obtain an integer $k_{1,7} \, (> k_{1,6})$ with $1!$ dividing $k_{1,7}$ and $2^2$ pairwise disjoint compact sets $K_{\al_0\al_1}^{(1,1,7)} \subset K_{\al_0\al_1}^{(1,1,6)}$, $\al_0 = 0, 1$, $\al_1 = 0, 1$, with nonempty interiors such that $1!$ divides $k_{1,7}$ and \\
[3pt]$f^{k_{1,7}}\big(K_{00}^{(1,1,7)}\big) \subset W_{1,1}^{(1)}, \, f^{k_{1,7}}\big(K_{01}^{(1,1,7)}\big) \subset W_{1,2}^{(1)},  f^{k_{1,7}}\big(K_{10}^{(1,1,7)}\big) \subset W_{1,3}^{(1)}, \, f^{k_{1,7}}\big(K_{11}^{(1,1,7)}\big) \subset W_{1,4}^{(1)}$. 

For all integers $1 \le r \le 2^4$, we let $\xi(1,r) = (2 \cdot {\it 1}^2 + 3 \cdot {\it 1} + 1 )+r = 6 + r$.  Similar arguments as above imply that there exist positive integers $(k_{1,6} <) \, k_{1,7} < k_{1,8} < k_{1,9} < \cdots < k_{1,21} < k_{1,22}$ and compact sets (with nonempty interiors) $(K_{\al_0\al_1}^{(1,1,6)} \supset) \, K_{\al_0\al_1}^{(1,1,7)} \supset K_{\al_0\al_1}^{(1,1,8)} \supset K_{\al_0\al_1}^{(1,1,9)} \supset \cdots \supset K_{\al_0\al_1}^{(1,1,21)} \supset K_{\al_0\al_1}^{(1,1,22)}$, $\al_0 = 0,1$, $0 \le i \le 1$, such that, for all $1 \le r \le 2^4$, if we relable the sets $K_{00}^{(1,1,\xi(1,r))}, K_{01}^{(1,1,\xi(1,r))}, K_{10}^{(1,1,\xi(1,r))}, K_{11}^{(1,1,\xi(1,r))}$ (in {\it lexicographical increasing order in} $\al_0\al_1$) as $K_{\xi(1,r),1}^{(1)}, K_{\xi(1,r),2}^{(1)}, K_{\xi(1,r),3}^{(1)}, K_{\xi(1,r),4}^{(1)}$, then, for all $1 \le r \le 2^4$, we have \\
[3pt]\centerline{$1!$ divides $k_{1,\xi(1,r)}$ and $f^{k_{1,\xi(1,r)}}\big(K_{\xi(1,r),m}^{(1)}\big) \subset W_{r,m}^{(1)}$ for all integers $1 \le m \le 4$.} \\
Consequently, for any open sets $W_1^{(1)}, W_2^{(1)}, W_3^{(1)}, W_4^{(1)}$ in $\{ U_1, U_2 \}$, there exists an integer $1 \le r(1) \le 2^4$, depending on $W_m^{(1)}, 1 \le m \le 4$, such that \\
[5pt]\centerline{$1!$ divides $k_{1,\xi(1,r(1))}$ and $f^{k_{1,\xi(1,r(1))}}\big(K_{\xi(1,r(1)),m}^{(1)}\big) \subset W_m^{(1)}$ for all $1 \le m \le 4$.}

At the second stage, 

let $n_2$ be any positive integer such that $n_2 > k_{1,22}$.  For each $\al_0 = 0,1$, $\al_1 = 0, 1$, let $K_{\al_0\al_1}^{(1,2,0)} = K_{\al_0\al_1}^{(1,1,22)}$.  Then, as indicated above, the set
\begin{multline*}
$$
U_2 \setminus \big(K_{00}^{(1,2,0)} \cup K_{01}^{(1,2,0)} \cup K_{10}^{(1,2,0)} \cup K_{11}^{(1,2,0)}\big) \, = 
U_2 \setminus \big(K_{00}^{(1,1,22)} \cup K_{01}^{(1,1,22)} \cup K_{10}^{(1,1,22)} \cup K_{11}^{(1,1,22)}\big) \\ 
\supset \, U_2 \setminus \big(K_{00}^{(1,1,1)} \cup K_{01}^{(1,1,1)} \cup K_{10}^{(1,1,1)} \cup K_{11}^{(1,1,1)}\big) \,\,\, \text{is nonempty}.
$$
\end{multline*}
Let $K_{00}^{(2,2,0)}, K_{01}^{(2,2,0)}, K_{10}^{(2,2,0)}, K_{11}^{(2,2,0)}$ be any $2^2$ pairwise disjoint compact sets (with nonempty interiors) in $U_2 \setminus \big( K_{00}^{(1,2,0)} \cup K_{01}^{(1,2,0)} \cup K_{10}^{(1,2,0)} \cup K_{11}^{(1,2,0)}\big)$  \big(note that the number $j$ in the superscript $(j,*,*)$ of $K_{\al_0\al_1}^{(j,*,*)}$ indicates that these sets lie in $U_j, j = 1,2$ and the number 2 in the superscript $(*,2,*)$ of $K_{\al_0\al_1}^{(*,2,*)}$ means that these sets are obtained at the second stage and so on\big). 

\indent In the following, let $\big(K_{00}^{(1,2,0)}\big)^2$ denote $K_{00}^{(1,2,0)}, K_{00}^{(1,2,0)}$ and let $\big(V_{n_2}(a_2), V_{n_2}(b_2)\big)^8$ denote 8 copies of $(V_{n_2}(a_2), V_{n_2}(b_2))$ and so on.  

By applying Lemma 3(1) to the following $2 \cdot 2^4$ sets with nonempty interiors \big(note the subscripts of $K_{...}^{(j,2,0)}$'s, $j= 1, 2$, and note that here we not only consider the $2^2$ just found compact sets $K_{00}^{(1,2,0)}, K_{01}^{(1,2,0)}, K_{10}^{(1,2,0)}, K_{11}^{(1,2,0)}$ from $U_1$, we also consider the $2^2$ newly chosen compact sets $K_{00}^{(2,2,0)}, K_{01}^{(2,2,0)}, K_{10}^{(2,2,0)}, K_{11}^{(2,2,0)}$ from $U_2$\big) \\
[5pt]\centerline{$\big(K_{\al_0\al_1}^{(1,2,0)}\big)^2, \, \al_i = 0, 1, \, 0 \le i \le 1; \, \big(K_{\al_0\al_1}^{(2,2,0)}\big)^2, \, \al_i = 0, 1, \, 0 \le i \le 1 \,\, : \,\, \big(V_{n_2}(a_2), V_{n_2}(b_2)\big)^8$,} \\
[3pt]\begin{footnotesize}
i.e., $(K_{00}^{(1,2,0)})^2, (K_{01}^{(1,2,0)})^2, (K_{10}^{(1,2,0)})^2, (K_{11}^{(1,2,0)})^2$; $(K_{00}^{(2,2,0)})^2, (K_{01}^{(2,2,0)}\big)^2, \big(K_{10}^{(2,2,0)})^2, (K_{11}^{(2,2,0)})^2$ : $(V_{n_2}(a_2), V_{n_2}(b_2))^8$,  
\end{footnotesize} 
we obtain an integer $k_{2,1} \, (> n_2 > k_{1,22})$ \big(the number 2 in the subscript $(2,*)$ of $k_{2,*}$ indicates that this number $k_{2,*}$ is obtained at the second stage\big) with $2!$ dividing $k_{2,1}$ and $2 \cdot 2^3$ pairwise disjoint compact sets $K_{\al_0\al_1\al_2}^{(j,2,1)} \subset K_{\al_0\al_1}^{(j,2,0)}, \, \al_i = 0, 1, \, 0 \le i \le 2, \, j = 1, 2$, with nonempty interiors such that \\
[10pt]\hspace*{.4in}each diameter is so smaller than $1/{2^5}$ that \\  
[4pt]\hspace*{1in}the open set $U_3 \setminus \big( \bigcup \big\{ K_{\al_0\al_1\al_2}^{(j,2,1)}: \, \al_i = 0, 1, \, 0 \le i \le 2, \, j = 1, 2 \big\} \big)$ is nonempty and \\
[6pt]\hspace*{.4in}$f^{k_{2,1}}\big(\bigcup_{\al_i = 0, 1, \, 0 \le i \le 2} \, (K_{\al_0\al_1\al_2}^{(1,2,1)} \cup K_{\al_0\al_1\al_2}^{(2,2,1)})\big) \subset W_{n_2}(\al_2)$, \\
[4pt]\hspace*{2in}where $W_{n_2}(0) = V_{n_2}(a_2)$ and $W_{n_2}(1) = V_{n_2}(b_2)$, and $2!$ divides $k_{2,1}$. \\
[12pt]\indent By applying Lemma 3(1) to the following $2 \cdot 2^4$ sets with nonempty interiors \\
[8pt]\centerline{$K_{\al_0\al_1\al_2}^{(1,2,1)}, \, \al_i = 0, 1, \, 0 \le i \le 2, \,\, K_{\al_0\al_1\al_2}^{(2,2,1)}, \, \al_i = 0, 1, \, 0 \le i \le 2 \,\, : \,\, \big(V_{n_2}(v_1)\big)^{16}$,} \\
[8pt]we obtain an integer $k_{2,2} \, (> k_{2,1})$ with $2!$ dividing $k_{2,2}$ and $2 \cdot 2^3$ pairwise disjoint compact sets $K_{\al_0\al_1\al_2}^{(j,2,2)} \subset K_{\al_0\al_1\al_2}^{(j,2,1)}, \, \al_i = 0, 1, \, 0 \le i \le 2, \, j = 1, 2$, with nonempty interiors such that 
$$
f^{k_{2,2}}\big(\cup \big\{ K_{\al_0\al_1\al_2}^{(j,2,2)}: \, \al_i = 0, 1, \, 0 \le i \le 2, \, j = 1, 2 \big\}\big) \, \subset V_{n_2}(v_1) \,\,\, \text{and} \,\,\, 2! \,\,\, \text{divides} \,\,\, k_{2,2}.
$$
\indent By applying Lemma 3(1) to the following $2 \cdot 2^4$ sets with nonempty interiors \\
[8pt]\centerline{$K_{\al_0\al_1\al_2}^{(1,2,2)}, \, \al_i = 0, 1, \, 0 \le i \le 2, \,\, K_{\al_0\al_1\al_2}^{(2,2,2)}, \, \al_i = 0, 1, \, 0 \le i \le 2 \,\, : \,\, \big(V_{n_2}(v_2)\big)^{16}$,} \\
[3pt]\begin{footnotesize}
i.e., $K_{000}^{(1,2,1)}, K_{001}^{(1,2,1)}, K_{010}^{(1,2,1)}, K_{011}^{(1,2,1)}, \cdots, K_{111}^{(1,2,1)}; \,\,K_{000}^{(2,2,1)}, K_{001}^{(2,2,1)}, K_{010}^{(2,2,1)}, K_{011}^{(2,2,1)}, \cdots, K_{111}^{(2,2,1)}$; \,\,$\big(V_{n_2}(v_2)\big)^{16}$,
\end{footnotesize} \\
[3pt]we obtain an integer $k_{2,3} \, (> k_{2,2})$ with $2!$ dividing $k_{2,3}$ and $2 \cdot 2^3$ pairwise disjoint compact sets $K_{\al_0\al_1\al_2}^{(j,2,3)} \subset K_{\al_0\al_1\al_2}^{(j,2,2)}, \, \al_i = 0, 1, \, 0 \le i \le 2, \, j = 1, 2$, with nonempty interiors such that 
$$
f^{k_{2,3}}\big(\cup \big\{ K_{\al_0\al_1\al_2}^{(j,2,3)}: \, \al_i = 0, 1, \, 0 \le i \le 2, \, j = 1, 2 \big\}\big) \, \subset V_{n_2}(v_2) \,\,\, \text{and} \,\,\, 2! \,\,\, \text{divides} \,\,\, k_{2,3}.
$$
\indent By applying Lemma 3(2) to the following $2 \cdot 2^4+1$ sets with nonempty interiors and $s = 0$ \\
[8pt]\centerline{$K_{\al_0\al_1\al_2}^{(j,2,3)}, \, \al_i = 0, 1, \, 0 \le i \le 2, \, j = 1, 2 \,\, : \,\, \big(G_{0,n_2}(x_1)\big)^{16}; \,\, G_{0,n_2}(x_1)$} \\
[8pt]we obtain an integer $k_{2,4} \, (> k_{2,3})$ and $2 \cdot 2^3$ pairwise disjoint compact sets $K_{\al_0\al_1\al_2}^{(j,2,4)} \subset K_{\al_0\al_1\al_2}^{(j,2,3)}$, $\al_0 = 0, 1$, $\al_1 = 0, 1$, $j = 1, 2$, with nonempty interiors such that 
$$
f^{k_{2,4}}\big(\cup \big\{ K_{\al_0\al_1\al_2}^{(j,2,4)}: \, \al_i = 0, 1, \, 0 \le i \le 2, \, j = 1, 2 \big\}\big) \, \subset G_{0,n_2}(x_1) \,\,\, \text{and} \,\,\, f^{k_{2,4}}(x_1) \in G_{0,n_2}(x_1).
$$
\indent By arguments similar to the previous paragraph and those steps from the third to the sixth at the first stage, we obtain 
\noindent
integers $(k_{2,4} <) \, k_{2,5} \, < k_{2,6} \, < k_{2,7} \, < \cdots < k_{2,14} \, < k_{2,15} \, $ and $2 \cdot 2^3$ pairwise disjoint compact sets ($K_{\al_0\al_1\al_2}^{(j,2,4)} \supset) \, K_{\al_0\al_1\al_2}^{(j,2,5)} \supset K_{\al_0\al_1\al_2}^{(j,2,6)} \supset K_{\al_0\al_1\al_2}^{(j,2,7)} \, \supset \cdots \supset K_{\al_0\al_1\al_2}^{(j,2,14)} \supset K_{\al_0\al_1\al_2}^{(j,2,15)}$, $\al_i = 0, 1$, $0 \le i \le 2$, $j = 1, 2$, with nonempty interiors such that, for $s = 0, 1, 2$, we have \\
[8pt]$f^{k_{2,4+s}}\big(\cup \big\{ f^s\big(K_{\al_0\al_1\al_2}^{(j,2,4+s)}\big): \, \al_i = 0, 1, \, 0 \le i \le 2, \, j = 1, 2 \big\}\big) \subset G_{s,n_2}(x_1) \, \text{and} \, f^{k_{2,4+s}}(x_1) \in G_{s,n_2}(x_1)$, \\
[5pt]$f^{k_{2,7+s}}\big(\cup \big\{ f^s\big(K_{\al_0\al_1\al_2}^{(j,2,7+s)}\big): \, \al_i = 0, 1, \, 0 \le i \le 2, \, j = 1, 2 \big\}\big) \subset \widehat G_{s,n_2}(x_1) \, \text{and} \, f^{k_{2,7+s}}(x_1) \in G_{s,n_2}(x_1)$, \\
[5pt]$f^{k_{2,10+s}}\big(\cup \big\{ f^s\big(K_{\al_0\al_1\al_2}^{(j,2,10+s)}\big): \, \al_i = 0, 1, \, 0 \le i \le 2, \, j = 1, 2 \big\}\big) \subset G_{s,n_2}(x_2) \, \text{and} \, f^{k_{2,10+s}}(x_2) \in G_{s,n_2}(x_2)$, \\
[5pt]$f^{k_{2,13+s}}\big(\cup \big\{ f^s\big(K_{\al_0\al_1\al_2}^{(j,2,13+s)}\big): \, \al_i = 0, 1, \, 0 \le i \le 2, \, j = 1, 2 \big\}\big) \subset \widehat G_{s,n_2}(x_2) \, \text{and} \, f^{k_{2,13+s}}(x_2) \in G_{s,n_2}(x_2)$. 

It is clear that there are $({\it 2}+1)^{2 \cdot 2^{{\it 2}+1}} = 3^{16}$ finite sequences of $2 \cdot 2^{{\it 2}+1}=16$ elements from the set $\{ U_1, U_2, U_3 \}$.  We label them as $W_{1,1}^{(2)}, W_{1,2}^{(2)}, W_{1,3}^{(2)}, \cdots, W_{1,15}^{(2)}, W_{1,16}^{(2)}$; $W_{2,1}^{(2)}, W_{2,2}^{(2)}, W_{2,3}^{(2)}, \cdots, W_{2,15}^{(2)}, W_{2,16}^{(2)}$; $W_{3,1}^{(2)}, W_{3,2}^{(2)}, W_{3,3}^{(2)}, \cdots, W_{3,15}^{(2)}, W_{3,16}^{(2)}; \,\,\, \cdots, \,\, W_{3^{16},1}^{(2)}, W_{3^{16},2}^{(2)}, W_{3^{16},3}^{(2)}$, $\cdots$; $W_{3^{16},15}^{(2)}, W_{3^{16},16}^{(2)}$. 

Now, for all integers $1 \le r \le 3^{16}$, we let $\xi(2,r) = (2\cdot {\it 2}^2 + 3 \cdot {\it 2} + 1) + r = 15+r$.  By firstly applying Lemma 3(1) to the following $2 \cdot 2^4$ sets with nonempty interiors \\
[8pt]\centerline{$K_{\al_0\al_1\al_2}^{(1,2,15)}, \, \al_i = 0, 1, \, 0 \le i \le 2, \,\, K_{\al_0\al_1\al_2}^{(2,2,15)}, \, \al_i = 0, 1, \, 0 \le i \le 2 \,\, : \,\, W_{1,m}^{(2)}, 1 \le m \le 16$,} \\
[8pt]and then by arguments similar to those steps from the $7^{th}$ to the $22^{nd}$ at the first stage, we obtain integers $(k_{2,15} <) \, k_{2,16} < k_{2,17} < k_{2,18} < \cdots < k_{2,15+3^{16}}$ and compact sets (with nonempty interiors) $(K_{\al_0\al_1\al_2}^{(j,2,15)} \supset) K_{\al_0\al_1\al_2}^{(j,2,16)} \supset K_{\al_0\al_1\al_2}^{(j,2,17)} \supset K_{\al_0\al_1\al_2}^{(j,2,18)} \supset \cdots \supset K_{\al_0\al_1\al_2}^{(j,2,15+3^{16})}$, $\al_i = 0,1, 0 \le i \le 2, j = 1, 2$, such that, for all $1 \le r \le 3^{16}$, if we relabel the sets $K_{\al_0\al_1\al_2}^{(j,2,\xi(2,r))}, \al_i = 0, 1, 0 \le i \le 2$, in {\it lexicographical increasing order in} $\al_0\al_1\al_2$ successively for $j = 1, 2$, as $K_{\xi(2,r),1}^{(2)}$, $K_{\xi(2,r),2}^{(2)}$, $\cdots$, $K_{\xi(2,r),16}^{(2)}$, then, for any $1 \le r \le 3^{16}$, we have \\
[5pt]\centerline{$2!$ divides $k_{2,\xi(2,r)}$ and $f^{k_{2,\xi(2,r)}}\big(K_{\xi(2,r),m}^{(2)}\big) \subset W_{r,m}^{(2)}$ for all integers $1 \le m \le 2^4$.}\\
[5pt]Consequently, for any open sets $W_1^{(2)}, W_2^{(2)}, W_3^{(2)}, \cdots, W_{16}^{(2)}$ in $\{ U_1, U_2, U_3 \}$, there exists an integer $1 \le r(2) \le 3^{16}$, depending on $W_m^{(2)}, 1 \le m \le 2^4$, such that \\
[5pt]\centerline{$2!$ divides $k_{2,\xi(2,r(2))}$ and $f^{k_{2,\xi(2,r(2))}}\big(K_{\xi(2,r(2)),m}^{(2)}\big) \subset W_m^{(2)}$ for all $1 \le m \le 2^4$.}

For simplicity, for all integers $\ell \ge 1$ and $1 \le r \le (\ell +1)^{\ell \cdot 2^{\ell +1}}$, we let $\eta(\ell) = (\ell +1)^{\ell \cdot 2^{\ell +1}}$ and $\xi(\ell,r) = (2\ell^2+3\ell+1) + r$ $\le (2\ell^2+3\ell+1) + \eta(\ell) = \xi(\ell,\eta(\ell))$.  Inductively, assume that we have completed the construction of the $(\ell-1)^{st}$ stage for some $\ell \ge 2$ and obtained integers $(n_{\ell-1} <) k_{\ell-1,1} < k_{\ell-1,2} < k_{\ell-1,3} < \cdots < k_{\ell-1,\xi(\ell-1,\eta(\ell-1))}$ and corresponding pairwise disjoint compact sets (with nonempty interiors) $\big($when $\ell \ge 3, K_{\al_0\al_1\cdots \al_{\ell-2}}^{(j, \, \ell-2, \, \xi(\ell-2, \eta(\ell-2)))} \supset$\big) $K_{\al_0\al_1 \cdots \al_{\ell-1}}^{(j,\, \ell-1,\, 1)} \supset K_{\al_0\al_1 \cdots \al_{\ell-1}}^{(j,\, \ell-1,\, 2)} \supset K_{\al_0\al_1 \cdots \al_{\ell-1}}^{(j,\, \ell-1,\, 3)} \supset \cdots \subset K_{\al_0\al_1 \cdots \al_{\ell-1}}^{(j,\, \ell-1,\, \xi(\ell-1, \eta(\ell-1)))}$, $\al_i = 0, 1$, $0 \le i \le \ell-1$, $1 \le j \le \ell-1$, which satisfy all relevant properties and the diameter of each such compact set is so smaller than $1/{2^{2\ell-1}}$ that the open set $U_{\ell} \setminus \ \big( \bigcup \big\{ K_{\al_0\al_1\al_2\cdots \al_{\ell-1}}^{(j, \,\ell-1, \, 1)}: \, \al_i = 0, 1, \, 0 \le i \le \ell-1, \, 1 \le j \le \ell-1 \big\} \big)$ is nonempty.  

For $\ell \ge 2$, at the $\ell^{th}$ stage, there are $(\ell +1)^{\ell \cdot 2^{\ell +1}}$ $=\eta(\ell)$ distinct sequences of elements in $\{U_1, U_2, \cdots, U_{\ell+1}\}$ of length $\ell \cdot 2^{\ell +1}$.  Label them as 
$W_{1,1}^{(\ell)}, W_{1,2}^{(\ell)}, \cdots, W_{1,\ell \cdot 2^{\ell +1}}^{(\ell)}$; 
$W_{2,1}^{(\ell)}, W_{2,2}^{(\ell)}, \cdots, W_{2,\ell \cdot 2^{\ell +1}}^{(\ell)}$; 
$\cdots, W_{\eta(\ell),1}^{(\ell)}$, $W_{\eta(\ell),2}^{(\ell)}$, $\cdots$, $W_{\eta(\ell),\ell \cdot 2^{\ell +1}}^{(\ell)}$.  Let $n_\ell$ be any positive integer such that $n_\ell > k_{\ell-1, \xi(\ell-1, \eta(\ell-1))}$.  For $1 \le j \le \ell-1, \al_i = 0, 1, 0 \le i \le \ell-1$, let $K_{\al_0\al_1\cdots \al_{\ell-1}}^{(j, \, \ell, \, 0)} = K_{\al_0\al_1\cdots \al_{\ell-1}}^{(j, \, \ell-1, \, \xi(\ell-1, \eta(\ell-1)))}$.  On the other hand, let $K_{\al_0\al_1\cdots \al_{\ell-1}}^{(\ell, \, \ell, \, 0)}$, $\al_i = 0, 1, 0 \le i \le \ell-1$, be any pairwise disjoint compact sets (with nonempty interiors) in $U_{\ell} \setminus \ \big( \bigcup \big\{ K_{\al_0\al_1\al_2\cdots \al_{\ell-1}}^{(j, \,\ell-1, \, 1)}: \, \al_i = 0, 1, \, 0 \le i \le \ell-1, \, 1 \le j \le \ell-1 \big\} \big)$.

By proceeding as the previous stages starting with the pairwise disjoint sets $K_{\al_0\al_1\cdots \al_{\ell-1}}^{(j, \, \ell, \, 0)}, \al_i = 0, 1, 0 \le i \le \ell-1, 1 \le j \le \ell$, we obtain $\xi(\ell,\eta(\ell))$ positive integers $(k_{\ell-1,\, \xi(\ell-1,\eta(\ell-1))} < n_\ell <) \, k_{\ell,1} < k_{\ell,2} < \cdots < k_{\ell,\,2\ell^2+3\ell+1} < k_{\ell,\,(2\ell^2+3\ell+1)+1} < k_{\ell,\,(2\ell^2+3\ell+1)+2} < \cdots < k_{\ell,\,\xi(\ell,\eta(\ell))}$ with ${\ell}!$ dividing $k_{\ell,i}$ for all $1 \le i \le \ell+1$ and all $(2\ell^2+3\ell+1)+1 \le i \le (2\ell^2+3\ell+1)+ (\ell+1)^{\ell \cdot 2^{\ell+1}} = \xi(\ell,\eta(\ell))$ and $\xi(\ell,\eta(\ell))$ {\it corresponding} pairwise disjoint compact sets (with nonempty interiors) $K_{\al_0\al_1\al_2\cdots \al_\ell}^{(j, \,\ell, \, r)}$, $\al_i = 0, 1, 0 \le i \le \ell$, $1 \le j \le \ell$, $1 \le r \le \xi(\ell,\eta(\ell))$ such that \\
[5pt]\hspace*{.2in}(when $\ell \ge 2$ and $\ell-1 \ge j \ge 1$, $(K_{\al_0\al_1\al_2\cdots \al_{\ell-1}}^{(j,\, \ell-1, \xi(\ell-1,\eta(\ell-1))} \supset) \,\, K_{\al_0\al_1\al_2\cdots \al_\ell}^{(j, \,\ell, \, 1)} \supset K_{\al_0\al_1\al_2\cdots \al_\ell}^{(j, \,\ell, \, 2)} \supset \cdots$ \\
[5pt]\hspace*{.8in}$\supset K_{\al_0\al_1\al_2\cdots \al_\ell}^{(j, \,\ell, \, 2\ell^2+3\ell+1)} \supset K_{\al_0\al_1\al_2\cdots \al_\ell}^{(j, \,\ell, \, (2\ell^2+3\ell+1)+1)} \supset K_{\al_0\al_1\al_2\cdots \al_\ell}^{(j, \,\ell, \, (2\ell^2+3\ell+1)+2)} \supset \cdots \supset K_{\al_0\al_1\al_2\cdots \al_\ell}^{(j, \,\ell, \, \xi(\ell,\eta(\ell))})$ \\
[5pt]and each diameter is so smaller than $1/{2^{2\ell+1}}$ that the open set \\
[5pt]\hspace*{.5in}$U_{\ell+1} \setminus \ \big( \bigcup \big\{ K_{\al_0\al_1\al_2\cdots \al_\ell}^{(j, \,\ell, \, 1)}: \, \al_i = 0, 1, \, 0 \le i \le \ell, \, 1 \le j \le \ell \big\} \big)$ is nonempty $\hfill(\star)$ \\
[6pt]and the following hold:  

\noindent
$f^{k_{\ell,1}}\big(\bigcup_{j=1}^\ell \big\{ \big(K_{\al_0\al_1\al_2\cdots\al_\ell}^{(j,\,\ell,\, \xi(\ell,\eta(\ell)))}\big): \, \al_i = 0, 1, \, 0 \le i \le \ell \big\}\big) \subset W_{n_\ell}(\al_\ell)$,\\ 
[5pt]\hspace*{2.5in}where $W_{n_\ell}(0) = V_{n_\ell}(a_\ell)$ and $W_{n_\ell}(1) = V_{n_\ell}(b_\ell)$, and,\\
[5pt]for each $1 \le m \le \ell$, 

\centerline{$f^{k_{\ell,1+m}}\big(\bigcup_{j=1}^\ell \big\{ K_{\al_0\al_1\al_2\cdots\al_\ell}^{(j, \, \ell, \,\xi(\ell,\eta(\ell)))}: \, \al_i = 0, 1, \, 0 \le i \le \ell \big\}\big) \subset V_{n_\ell}(v_m)$ and,}

\noindent
for each $1 \le m \le \ell$,  \\
[5pt]\hspace*{.3in}$f^{k_{\ell,(1+\ell)+(2m-2)(\ell+1)+1+s}}\big(\bigcup_{j=1}^\ell \big\{ f^s\big(K_{\al_0\al_1\al_2\cdots\al_\ell}^{(j, \, \ell,\, \xi(\ell,\eta(\ell)))}\big): \, \al_i = 0, 1, \, 0 \le i \le \ell \big\}\big)$ \\ 
[5pt]\hspace*{.8in}$\subset G_{s,n_\ell}(x_m)$ and $f^{k_{\ell,(1+\ell)+(2m-2)(\ell+1)+1+s}}(x_m) \in G_{s,n_\ell}(x_m)$ for all $0 \le s \le \ell$, and\\
[5pt]\hspace*{.3in}$f^{k_{\ell,(1+\ell)+(2m-1)(\ell+1)+1+s}}\big(\bigcup_{j=1}^\ell \big\{f^s\big(K_{\al_0\al_1\al_2\cdots\al_\ell}^{(j, \, \ell, \,\xi(\ell,\eta(\ell)))}\big): \, \al_i = 0, 1, \, 0 \le i \le \ell \big\}\big)$ \\ 
[5pt]\hspace*{.8in}$\subset \widehat G_{s,n_\ell}(x_m)$ and $f^{k_{\ell,(1+\ell)+(2m-1)(\ell+1)+1+s}}(x_m) \in G_{s,n_\ell}(x_m)$ for all $0 \le s \le \ell$, and  

\noindent furthermore, 
\begin{quote}
\small
For each $1 \le r \le (\ell+1)^{\ell \cdot 2^{\ell+1}}$, by relabelling the sets $K_{\al_0\al_1\al_2\cdots \al_\ell}^{(j, \, \ell, \, \xi(\ell,r))}$, $\al_i = 0, 1$, $0 \le i \le \ell$, in {\it lexicographical increasing order in} $\al_0\al_1\al_2 \cdots \al_\ell$ successively for $1 \le j \le \ell$ as $K_{\xi(\ell,r),1}^{(\ell)}$, $K_{\xi(\ell,r),2}^{(\ell)}$, $\cdots$, $K_{\xi(\ell,r),\ell \cdot 2^{\ell+1}}^{(\ell)}$.  Then, for any open sets $W_1^{(\ell)}, W_2^{(\ell)}, \cdots, W_{\ell \cdot 2^{\ell+1}}^{(\ell)}$, in $\{ U_1, U_2, \cdots, U_{\ell+1}\}$, there is an integer $1 \le r(\ell) \le (\ell+1)^{\ell \cdot 2^{\ell+1}}$, depending on $\ell$ and $W_m^{(\ell)}, 1 \le m \le \ell \cdot 2^{\ell+1}$, such that $\ell!$ divides $k_{\ell,\xi(\ell,r(\ell))}$ and $f^{k_{\ell,\xi(\ell,r(\ell))}}\big(K_{\xi(\ell,r(\ell)),m}^{(\ell)}\big) \subset W_m^{(\ell)}$ for all $1 \le m \le \ell \cdot 2^{\ell+1}$.  $\hfill(\star\star)$
\end{quote}

\indent Let $\mathcal M = \{ \, k_{\ell,i}: 1 \le i \le \ell+1 \, \text{and} \, (2\ell^2+3\ell+1)+ 1 \le i \le (2\ell^2+3\ell+1)+ (\ell +1)^{\ell \cdot 2^{\ell+1}}, \ell \ge 1 \}$.  Then it follows from our construction that, for each integer $n \ge 1$, $n !$ (and so, $n$) divides all sufficiently large integers in $\mathcal M$.  Recall that $\xi(\ell,\eta(\ell)) = (2\ell^2+3\ell+1)+ (\ell+1)^{\ell \cdot 2^{\ell+1}}$.  

Let $\Sigma_2 = \{ \al : \al = \al_0\al_1\al_2 \cdots, \al_i = 0, 1, \, i \ge 0 \}$ and, for any $\al \in \Sigma_2$ and any $j \ge 1$, let \\
[6pt]\centerline{$K_{\al}^{(j)} = \bigcap_{\ell \, \ge \, j} \, K_{\al_0\al_1\al_2 \cdots \al_\ell}^{(j, \,\ell, \,\xi(\ell,\eta(\ell)))} \, \big(\subset U_j\big)$.} \\
[6pt]Then, for each $\al \in \Sigma_2$, since the sequence $< K_{\al_0\al_1\al_2\cdots \al_\ell}^{(j,\, \ell, \,\xi(\ell,\eta(\ell)))}>_{\ell \ge j}$ of compact sets is nested and their diameters shrink to zero, each $K_{\al}^{(j)}$ consists of exactly one point, say $K_{\al}^{(j)} = \{ x_{\al}^{(j)} \}$.  For each $j \ge 1$, let $\mathcal S^{(j)} = \{ x_{\al}^{(j)} : \al \in \Sigma_2 \} = \bigcap_{\ell \ge j} \big\{ \bigcup \, \{ K_{\al_0\al_1\al_2 \cdots \al_\ell}^{(j, \, \ell, \, \xi(\ell,\eta(\ell)))}: \al_i = 0, 1, 0 \le i \le \ell \}\big\} \, \big(\subset U_j\big)$.  Then $\mathcal S^{(j)} \subset U_j$ for each $j \ge 1$ and since $\{U_1, U_2, U_3, \cdots \}$ is a countable open base for $\mathcal X$, the set $\mathbb S = \bigcup_{j=1}^\infty \, \mathcal S^{(j)}$ is dense in $\mathcal X$.  

It is easy to see that $\mathcal S^{(j)}$ is a Cantor set.  Indeed, since $\mathcal S^{(j)}$ is the intersection of a nested sequence of compact sets, it is compact.  On the other hand, let $x_\al^{(j)}$ be a point in $\mathcal S^{(j)}$ and let $U$ be an open neighborhood of $x_\al^{(j)}$.  Since $K_\al^{(j)} = \{ x_\al^{(j)} \}$, there exists a positive integer $\hat \ell$ such that $x_\al^{(j)} \in K_{\al_0\al_1\al_2\cdots\al_{\hat \ell}}^{(j,\, \hat \ell,\, \xi(\hat \ell, \eta(\hat \ell)))} \subset U$.  Since $K_{\al_0\al_1\al_2\cdots\al_{\hat \ell}}^{(j,\, \hat \ell,\, \xi(\hat \ell, \eta(\hat \ell)))}$ contains the points $x_{\al_0\al_1\al_2\cdots\al_{\hat \ell}\be_{\hat \ell+1}\be_{\hat \ell+2}\cdots}^{(j)}$ for all $\be_{\hat \ell+i} = 0, 1, i \ge 1$, $U$ contains uncountably many points of $\mathcal S^{(j)}$.  So, $\mathcal S^{(j)}$ is perfect.  Furthermore, let $x_{\al}^{(j)}$ and $x_{\be}^{(j)}$ be any two distinct points in $\mathcal S^{(j)}$ with $\al \ne \be$.  As before, there is a positive integer $\check \ell$ such that $x_\al^{(j)} \in K_{\al_0\al_1\al_2\cdots\al_{\check \ell}}^{(j,\, \check \ell,\, \xi(\check \ell,\eta(\check \ell)))}$ and $x_\be^{(j)} \in K_{\be_0\be_1\be_2\cdots\be_{\check \ell}}^{(j,\, \check \ell,\, \xi(\check \ell,\eta(\check \ell)))}$.  Since the sets $K_{\ga_0\ga_1\ga_2\cdots\ga_{\check \ell}}^{(j,\, \check \ell,\, \xi(\check \ell,\eta(\check \ell)))}, \ga_i = 0, 1, 0 \le i \le \check \ell, 1 \le i \le j$, are pairwise disjoint, there exists disjoint open sets $U$ and $V$ such that 
$x_\al^{(j)} \in K_{\al_0\al_1\al_2\cdots\al_{\check \ell}}^{(j,\, \check \ell,\, \xi(\check \ell,\eta(\check \ell)))} \subset U$ and $x_\be^{(j)} \in \big(\bigcup \{K_{\ga_0\ga_1\ga_2\cdots\ga_{\check \ell}}^{(j,\, \check \ell,\, \xi(\check \ell,\eta(\check \ell)))}: \ga_i = 0, 1, 0 \le i \le \check \ell \}\big) \setminus K_{\al_0\al_1\al_2\cdots\al_{\check \ell}}^{(j,\, \check \ell,\, \xi(\check \ell,\eta(\check \ell)))} \subset V$.  Therefore, $\mathcal S^{(j)}$ is totally disconnected.  Since we have shown that $\mathcal S^{(j)}$ is compact, perfect and totally disconnected, $\mathcal S^{(j)}$ is a Cantor set and it follows from the above ($\star$) that the Cantor sets $\mathcal S^{(1)}, \mathcal S^{(2)}, \mathcal S^{(3)}, \cdots$ are pairwise disjoint.

Let $n \ge 1$ be any integer.  Let $a_1, a_2, \cdots, a_n$ be $n$ distinct points in $\cup_{j=1}^\infty \mathcal S^{(j)}$ and let $W_1', W_2', W_3'$, $\cdots$, $W_n'$ be any $n$ nonempty open sets in $\mathcal X$.  For each integer $1 \le m \le n$, let $p_m$ be a positive integer such that $U_{p_{m}} \subset W_m'$ (the integers $p_m$'s need not be distinct).  Then, since, for any $j \ge 1$, any $\al = \al_0\al_1\al_2 \cdots \in \Sigma_2$, and any $\ell \ge j$, we have the nested sequence \\
[5pt]\noindent \,\,\, $K_{\al_0\al_1\al_2\cdots \al_{\ell}}^{(j, \, \ell, \, 1)} \supset K_{\al_0\al_1\al_2\cdots \al_{\ell}}^{(j, \, \ell, \, 2)} \supset \cdots \supset K_{\al_0\al_1\al_2\cdots \al_{\ell}}^{(j, \, \ell, \, 2\ell^2+3\ell+1)} \supset \cdots \supset K_{\al_0\al_1\al_2\cdots \al_{\ell}}^{(j, \, \ell, \, \xi(\ell,\eta(\ell)))}$ 
\begin{equation*}
\supset K_{\al_0\al_1\al_2\cdots \al_{\ell+1}}^{(j, \, \ell+1, \, 1)} \supset K_{\al_0\al_1\al_2\cdots \al_{\ell+1}}^{(j, \, \ell+1, \, 2)} \supset \cdots \supset K_{\al_0\al_1\al_2\cdots \al_{\ell+1}}^{(j, \, \ell+1, \, \xi(\ell+1,\eta(\ell+1)))} \supset \cdots \tag{$\star\star\star$} \\
\end{equation*}
which satisfies that $\bigcap_{\ell \ge j} K_{\al_0\al_1\al_2\cdots \al_{\ell}}^{(j, \, \ell, \, \xi(\ell,\eta(\ell)))} = \{ x_\al^{(j)} \}$, a set consists of exactly one point, there exists an integer $\hat \ell > \max \{ p_1, p_2, \cdots, p_n \}$ such that $\{ a_1, a_2, \cdots, a_n \} \subset \bigcup_{j=1}^{\hat \ell} \mathcal S^{(j)}$ and, for each $\al_i = 0, 1$, $0 \le i \le \hat \ell$, and each $1 \le j \le \hat \ell$, the set $K_{\al_0\al_1\al_2 \cdots \al_{\hat \ell}}^{(j,\, \hat \ell, \, \xi(\hat \ell,\eta(\hat \ell)))}$ contains {\it at most} one point of $\{ a_1, a_2, \cdots, a_n \}$.  So, for each $\ell > \hat \ell$, it follows from ($\star\star\star$) that, for each $\al_i = 0, 1$, $0 \le i \le \ell$, $1 \le j \le \ell$, and each $1 \le r \le (\ell+1)^{\ell \cdot 2^{\ell+1}}$, the set $K_{\al_0\al_1\al_2 \cdots \al_{\hat \ell}\cdots \al_{\ell}}^{(j,\, \ell, \, \xi(\ell,r))} \big(\subset K_{\al_0\al_1\al_2 \cdots \al_{\hat \ell}}^{(j,\, \hat \ell, \, \xi(\hat \ell,\eta(\hat \ell)))}\big)$ contains {\it at most} one point of $\{ a_1, a_2, \cdots, a_n \}$ and 
$\{a_1, a_2, \cdots, a_n\} \subset \cup_{j=1}^{\hat \ell} \mathcal S^{(j)} \subset \cup_{j=1}^{\hat \ell} \big\{K_{\al_0\al_1\al_2 \cdots \al_{\hat \ell}\cdots \al_{\ell}}^{(j,\, \ell, \, \xi(\ell,r))}: \al_i = 0, 1, 0 \le i \le \ell \big\}$.

For each $\ell > \hat \ell$, all $1 \le j \le \ell$ and all $1 \le r \le (\ell+1)^{\ell \cdot 2^{\ell+1}}$, we relabel the sets $K_{\al_0\al_1 \cdots \al_{\ell}}^{(j,\ell,\xi({\ell,r}))}, \al_i = 0,1, 0 \le i \le \ell$, in {\it lexicographical increasing order in} $\al_0\al_1\al_2\cdots \al_\ell$ successively for $1 \le j \le \ell$,  as $K_{\xi({\ell,r}),1}^{(\ell)}$, $K_{\xi({\ell,r}),2}^{(\ell)}$, $\cdots$, $K_{\xi({\ell,r}),\ell \cdot 2^{\ell+1}}^{(\ell)}$.  Then, for each $1 \le m \le n$, there exists an integer $1 \le t({\ell,m}) \le \ell \cdot 2^{\ell+1}$, such that \\
[6pt]\centerline{$a_{m} \in K_{\xi({\ell, \eta(\ell)}),t(\ell,m)}^{(\ell)} \big(\subset K_{\xi({\ell,r}),t({\ell,m})}^{(\ell)}$ for all $1 \le r \le \eta(\ell) = (\ell+1)^{\ell \cdot 2^{\ell+1}}$\big)} \\
[6pt]and these positive integers $t({\ell,1}), t({\ell,2}), t({\ell,3}), \cdots, t({\ell,n})$ are all distinct.

For each $\ell > \hat \ell$, let $W_{1}^{(\ell)}$, $W_{2}^{(\ell)}$, $\cdots$, $W_{\ell \cdot 2^{\ell+1}}^{(\ell)}$ be open sets in $\{U_1, U_2, \cdots, U_{\ell+1} \}$ such that 
$$
W_{t(\ell,m)}^{(\ell)} = U_{p_{m}} \subset W_m' \,\,\, \text{for all integers} \,\,\, 1 \le m \le n \,\,\, \text{and} \,\,\, W_{i}^{(\ell)} = U_1 \,\,\, \text{otherwise}.                  
$$
Then it follows from the above ($\star\star$) that there is an integer $1 \le r(\ell) \le (\ell+1)^{\ell \cdot 2^{\ell+1}}$ such that $\ell!$ divides $k_{\ell,\xi(\ell,r(\ell))}$ and $f^{k_{\ell,\xi(\ell,r(\ell))}} \big(K_{\xi(\ell,r(\ell)),i}^{(\ell)}\big) \subset W_{i}^{(\ell)}$ for all $1 \le i \le \ell \cdot 2^{\ell+1}$.  In particular, we obtain that $f^{k_{\ell,\xi(\ell,r(\ell))}}\big(K_{\xi(\ell,r(\ell)),t(\ell,m)}^{(\ell)}\big) \subset W_{t(\ell,m)}^{(\ell)} = U_{p_m} \subset W_m'$ for all $1 \le m \le n$.  Consequently, for each $1 \le m \le n$, we have $f^{k_{\ell,\xi(\ell,r(\ell))}}(a_{m}) \in f^{k_{\ell,\xi(\ell,r(\ell))}}\big(K_{\xi(\ell,\eta(\ell)),t(\ell,m)}^{(\ell)}\big) \subset f^{k_{\ell,\xi(\ell,r(\ell))}}\big(K_{\xi(\ell,r(\ell)),t(\ell,m)}^{(\ell)}\big) \subset W_{t(\ell,m)}^{(\ell)} = U_{p_{m}} \subset W_m'$. Since $k_{\ell,\xi(\ell,r(\ell))} \in \mathcal M$ for each $\ell > \hat \ell$ and $\ell!$ divides all sufficiently large numbers in $\mathcal M$, and since $W_1', W_2', \cdots, W_n'$ are arbitrary nonempty open sets in $\mathcal X$,  this confirms Part (1).   

As for Part (2), let $c$ and $d$ be any two distinct points in $\mathbb S = \bigcup_j^\infty \mathcal S^{(j)}$ and let $<(a_k, b_k)>_{k \ge 1}$ be a sequence of points in $\mathcal X \times \mathcal X$ such that $\lim_{k \to \infty} \rho(a_k, b_k) = \be$.  Then, it follows from Part (1) that the set $\{\big(f^m(c), f^m(d)\big): m \in \mathcal M \}$ is dense in $\mathcal X \times \mathcal X$.  Let $W_1$ and $W_2$ be nonempty open sets in $\mathcal X$ and let $s \ge 0$ be an integer.  Then $(f \times f)^{-s}(W_1 \times W_2)$ is a nonempty open set in $\mathcal X \times \mathcal X$ and so contains $(f^{m'}(c), f^{m'}(d))$ for some $m' \in \mathcal M$.  Consequently, for any integer $s \ge 0$, the set $\{\big(f^{m+s}(c), f^{m+s}(d)\big): m \in \mathcal M \}$ is dense in $\mathcal X \times \mathcal X$.  In particular, there exists a strictly increasing sequence $< m_i >_{i \ge 1}$ of positive integers in $\mathcal M$ such that $\lim_{i \to \infty} \big(f^{m_i+s}(c), f^{m_i+s}(d)\big) = (c, c)$.  Furthermore, for each $k \ge 1$, there exists a strictly increasing sequence $< m_{k,i} >_{i \ge 1}$ of positive integers in $\mathcal M$ such that $\lim_{i \to \infty} \big(f^{m_{k,i}+s}(c), f^{m_{k,i}+s}(d)\big) = (a_k, b_k)$.  By taking the diagonal sequence, we obtain that $\lim_{k \to \infty} \big(f^{m_{k,k}+s}(c), f^{m_{k,k}+s}(d)\big) = \lim_{k \to \infty} (a_k, b_k) = \beta$.  Consequently, we obtain that
\begin{equation*}
\liminf_{\substack{m \to \infty \\ m \in \mathcal M}} \rho\big(f^m\big(f^s(c)\big), f^m\big(f^s(d)\big)\big) = 0 \,\,\, \text{and} \,\,\, \limsup_{\substack{m \to \infty \\ m \in \mathcal M}} \rho\big(f^m\big(f^s(c)\big), f^m\big(f^s(d)\big)\big) = \be.  \tag{$\dagger$}
\end{equation*}
\indent Since $n$ divides all sufficiently large integers in $\mathcal M$, we obtain that Part (2) holds as an easy consequence of Part (1).  The fact ($\dagger$) will be used in proving Part (5) later.

For all positive integers $r$ and $m$, since \\   
[8pt]\centerline{$\mathcal S^{(1)} \cup \mathcal S^{(2)} \cup \cdots \cup \mathcal S^{(r)} \subset \bigcup_{j=1}^\ell \big\{ K_{\al_0\al_1\al_2\cdots\al_\ell}^{(j, \ell,\xi(\ell,\eta(\ell)))}: \al_i = 0, 1, \, 0 \le i \le \ell \big\} \, \text{for each integer} \, \ell \ge r$,} \\
[10pt]and since, for each integer $\ell \ge \max \{ r, m \}$, \\
[8pt]\hspace*{.3in}$k_{\ell,1+m} \in \mathcal M$ and $f^{k_{\ell,1+m}}\big(\bigcup_{j=1}^\ell \big\{K_{\al_0\al_1\al_2\cdots\al_\ell}^{(j,\ell,\xi(\ell,\eta(\ell)))}: \al_i = 0, 1, \, 0 \le i \le \ell \big\}\big) \subset \, V_{n_\ell}(v_m)$,\\
[8pt]and since $V_{n_\ell}(v_m)$ shrinks to the point $v_m$ as $n_\ell$ tends to $\infty$, we obtain that, for each integer $r \ge 1$, 
the union $\mathcal S^{(1)} \cup \mathcal S^{(2)} \cup \cdots \cup \mathcal S^{(r)}$ is $\mathcal M$-synchronouly proximal to $v_m$ for each $m \ge 1$.  This establishes Part (3). \\
[10pt]\indent For any integers $r \ge 1$, $m \ge 1$, $s \ge 0$ and all $\ell \ge \max \{ r, m, s \}$, it follows from the definitions of 
$k_{\ell,1+\ell+(2m-2)(\ell+1)+1+s}$ \, and \, $k_{\ell,1+\ell+(2m-1)(\ell+1)+1+s}$ and (respectively) the corresponding compact sets (with nonempty interiors) \\
[8pt]\hspace*{.2in}$K_{\al_0\al_1\al_2\cdots\al_\ell}^{(j,\ell,1+\ell+(2m-2)(\ell+1)+1+s)} \, \big(\supset K_{\al_0\al_1\al_2\cdots\al_\ell}^{(j,\ell,\xi(\ell,\eta(\ell)))}$, $\al_i = 0, 1$, $0 \le i \le \ell\big)$, $1 \le j \le \ell$ and \\ 
[6pt]\hspace*{.2in}$K_{\al_0\al_1\al_2\cdots\al_\ell}^{(j,\ell,1+\ell+(2m-1)(\ell+1)+1+s)} \, \big(\supset K_{\al_0\al_1\al_2\cdots\al_\ell}^{(j,\ell,\xi(\ell,\eta(\ell)))}$, $\al_i = 0, 1$, $0 \le i \le \ell\big)$, $1 \le j \le \ell$ \\
[8pt]that \\
[8pt]\hspace*{.675in}$f^{k_{\ell,1+\ell+(2m-2)(\ell+1)+1+s}}\big(f^s(\mathcal S^{(1)} \cup \mathcal S^{(2)} \cup \mathcal S^{(3)} \cup \cdots \cup \mathcal S^{(r)})\big)$ \\
[6pt]\hspace*{.5in}$\subset f^{k_{\ell,1+\ell+(2m-2)(\ell+1)+1+s}}\big(\bigcup_{j=1}^\ell \big\{ f^s(K_{\al_0\al_1\al_2\cdots\al_\ell}^{(j,\ell,\xi(\ell,\eta(\ell)))}: \, \al_i = 0, 1, \, 0 \le i \le \ell \big\}\big)$ \\ 
[6pt]\hspace*{.5in}$\subset f^{k_{\ell,1+\ell+(2m-2)(\ell+1)+1+s}}\big(\bigcup_{j=1}^\ell \big\{ f^s(K_{\al_0\al_1\al_2\cdots\al_\ell}^{(j,\ell,1+\ell+(2m-2)(\ell+1)+1+s)}): \, \al_i = 0, 1, \, 0 \le i \le \ell \big\}\big)$ \\ 
[6pt]\hspace*{.5in}$\subset G_{s,n_\ell}(x_m)$ and $f^{k_{\ell,1+\ell+(2m-2)(\ell+1)+1+s}}(x_m) \in G_{s,n_\ell}(x_m)$, \\
[8pt]and \\   
[8pt]\hspace*{.675in}$f^{k_{\ell,1+\ell+(2m-1)(\ell+1)+1+s}}\big(f^s(\mathcal S^{(1)} \cup \mathcal S^{(2)} \cup \mathcal S^{(3)} \cup \cdots \cup \mathcal S^{(r)})\big)$ \\
[6pt]\hspace*{.5in}$\subset f^{k_{\ell,1+\ell+(2m-1)(\ell+1)+1+s}}\big(\bigcup_{j=1}^\ell  \big\{ f^s(K_{\al_0\al_1\al_2\cdots\al_\ell}^{(j,\ell,\xi(\ell,\eta(\ell)))}: \, \al_i = 0, 1, \, 0 \le i \le \ell \big\}\big)$ \\ 
[6pt]\hspace*{.5in}$\subset f^{k_{\ell,1+\ell+(2m-1)(\ell+1)+1+s}}\big(\bigcup_{j=1}^\ell \big\{ f^s(K_{\al_0\al_1\al_2\cdots\al_\ell}^{(j,\ell,1+\ell+(2m-1)(\ell+1)+1+s)}): \, \al_i = 0, 1, \, 0 \le i \le \ell \big\}\big)$ \\ 
[6pt]\hspace*{.5in}$\subset \widehat G_{s,n_\ell}(x_m)$ and $f^{k_{\ell,1+\ell+(2m-1)(\ell+1)+1+s}}(x_m) \in G_{s,n_\ell}(x_m)$. \\
[8pt]\indent Therefore, we obtain that, for any point $x$ in $\mathcal S^{(1)} \cup \mathcal S^{(2)} \cup \mathcal S^{(3)} \cup \cdots \cup \mathcal S^{(r)}$, \\
[6pt]\hspace*{.5in}$f^{k_{\ell,1+\ell+(2m-2)(\ell+1)+1+s}}\big(\{f^s(x), x_m \}\big) \subset G_{s,n_\ell}(x_m)$ and \\
[6pt]\hspace*{.5in}$f^{k_{\ell,1+\ell+(2m-1)(\ell+1)+1+s}}\big(f^s(x)\big) \in \widehat G_{s,n_\ell}(x_m)$ and $f^{k_{\ell,1+\ell+(2m-1)(\ell+1)+1+s}}(x_m) \in G_{s,n_\ell}(x_m)$. \\
[8pt]\indent Since the diameter of the set $G_{s,n_\ell}(x_m)$ shrinks to zero as $n_\ell$ tends to $\infty$, we obtain that $\liminf_{n \to \infty} \rho\big(f^n(f^s(x)), f^n(x_m)\big) = 0$.  Furthermore, since, by the choices of $\widehat G_{s,n_\ell}(x_m)$'s, we have $\lim_{\ell \to \infty} dist\big(\widehat G_{s,n_\ell}(x_m), G_{s,\,n_\ell}(x_m)\big) \ge \beta/2$ and so, we obtain that $\limsup_{n \to \infty} \rho\big(f^n(f^s(x)), f^n(x_m)\big) \ge \beta/2$.  This verifies Part (4).

Now suppose $f$ has a fixed point $z$ in $\mathcal X$.  In the following, we shall use the above fact $(\dagger)$ with $s \ge 0$ to prove the existence of {\it invariant} $\delta$-scrambled sets of $f$ in $\mathcal X$ for some $\delta \ge 0$.

Let $v_1 = z$, a fixed point of $f$, $v_2 = u$, a transitive point of $f$, and, for each $n \ge 3$, let $v_n = f^{n-2}(u)$ (redundant).  Let $\mathcal S^{(1)}$, $\mathcal S^{(2)}$, $\mathcal S^{(3)}$, $\cdots$, be the resulting Cantor sets of totally transitive points of $f$, $\mathcal M$ the resulting infinite set of positive integers, $\mathbb S = \bigcup_{j = 1}^\infty \mathcal S^{(j)}$, the resulting $\beta$-scrambled sets (with respect to $\mathcal M$) of $f$ obtained by arguing as above with respect to the chosen set $\{ v_1, v_2, v_3, \cdots \}$.

Let $c$ and $d$ be any (may be identical) points in $\mathbb S$ and let $t \ge s \ge 0$ be any integers. 

Since, for each $\ell \ge 1$, the union $\bigcup_{j=1}^\ell \mathcal S^{(j)}$ is $\mathcal M$-synchronously proximal to $v_1 = z$, so is the set $\{ c, d \}$.  Let $< m_i>_{i \ge 1}$ be a strictly increasing sequence of positive integers in $\mathcal M$ such that $\lim_{m_i \to \infty} f^{m_i}(c) = v_1 = z = \lim_{m_i \to \infty} f^{m_i}(d)$.  Consequently, we have 
$0 \le \liminf_{\substack{m \to \infty \\ m \in \mathcal M}} \rho\big(f^m(f^t(c))$, $f^m(f^s(d))\big) \le \liminf_{i \to \infty} \rho\big(f^t(f^{m_i}(c)), f^s(f^{m_i}(d))\big) = \rho(z, z) = 0$.  This implies that
\begin{equation*}
\liminf_{\substack{m \to \infty \\ m \in \mathcal M}} \rho\big(f^m(f^t(c)), f^m(f^s(d))\big) = 0. \tag{$\ddagger$}
\end{equation*}
\indent Suppose $t > s \ge 0$.  Since $v_2 = u$ is a transitive point of $f$, we obtain that, for any integers $0 \le n_1 < n_2$, $f^{n_1}(u) \ne f^{n_2}(u)$.  Since, for each $\ell \ge 1$, the union $\bigcup_{j=1}^\ell \mathcal S^{(j)}$ is $\mathcal M$-synchronously proximal to $v_2 = u$, there exist a strictly increasing sequence $< m_k' >_{k \ge 1}$ of positive integers in $\mathcal M$ such that, for each $k \ge 1$, both the points $f^{m_k'}(c)$ and $f^{m_k'}(d)$ are {\it so close} to $u$ that 
\begin{multline*}
$$
f^{m_k'}\big(f^t(c)\big) = f^t\big(f^{m_k'}(c)\big) \approx f^t(u), \, f^{m_k'}\big(f^s(d)\big) = f^s\big(f^{m_k'}(d)\big) \approx f^s(u) \,\,\, \text{and} \\ 
\rho\big(f^{m_k'}(f^t(c)), \, f^{m_k'}(f^s(d))\big) \approx \rho\big(f^t(u), f^s(u)\big) > 0. 
$$
\end{multline*}
Therefore, we obtain that $\limsup_{\substack{m \to \infty \\ m \in \mathcal M}} \rho\big(f^m(f^t(c)), f^m(f^s(d))\big) \ge \lim_{k \to \infty} \rho\big(f^{m_k'}(f^t(c)), \, f^{m_k'}(f^s(d))\big)$ $> 0$ when $t > s \ge 0$.  On the other hand, for $t = s \ge 0$, $c \ne d$ and $\{ c, d \} \subset \mathbb S$, it follows from ($\dagger$) that $\limsup_{\substack{m \to \infty \\ m \in \mathcal M}} \rho\big(f^m(f^s(c)), f^m(f^s(d))\big) = \beta > 0$. 

These, combined with the above $(\ddagger)$, imply that the set $\widehat {\mathbb S} = \bigcup_{i=0}^\infty \, f^i(\mathbb S)$ is a dense {\it invariant} $0$-scrambled set (with respect to $\mathcal M$) of totally transitive points of $f$ in $\mathcal X$.  Since, for each $n \ge 1$, $n$ divides all sufficiently large positive integers in $\mathcal M$, we obtain that the set $\widehat {\mathbb S}$ is an {\it invariant} $0$-scrambled set of $f^n$ for each $n \ge 1$.

Finally, suppose $f$ has a fixed point $z$ and $\delta = \inf \big\{ \sup\{ \rho(f^n(x), x): x \in \mathcal X \}: n \ge 1 \big\}$ is $> 0$.  Then $f$ is not uniformly rigid.  It follows from {\bf\cite{fo}} that $f$ has an {\it invariant} $\varepsilon$-scrambled set for some $\varepsilon > 0$.  Here, we show that there exists a dense {\it invariant} set of $f$ in $\mathcal X$ which is a $\delta$-scrambled set of $f^n$ for all $n \ge 1$.

For any positive integers $n$ and $k$, let $\delta_n = \sup\{ \rho(f^n(x), x): x \in \mathcal X \}$ and let $c_{n,\,k}$ be a point in $\mathcal X$ such that $\rho\big(f^n(c_{n,\,k}), c_{n,\,k} \big) > \delta_n - 1/k$.  Since $f$ is a continuous weakly mixing map on $\mathcal X$, it follows from the above that $f$ has totally transitive points in $\mathcal X$.  Therefore, for each integer $s \ge 1$, $f^s$ is transitive on $\mathcal X$ and so, $f^s(\mathcal X)$ is dense in $\mathcal X$.  Consequently, for each $n \ge 1$, each $s \ge 1$ and each $k \ge 1$, there exists a point $c_{n,s,k}$ in $\mathcal X$ such that the point $f^s(c_{n,s,k})$ is so close to $c_{n,k}$ that 
\begin{equation*}
\rho\big(f^n(f^s(c_{n,s,k})), f^s(c_{n,s,k})\big) > \delta_n - 1/k. \tag{$\dagger\dagger$}   
\end{equation*}
\indent For each $n \ge 1$ and each $k \ge 1$, let $c_{n,0,k} = c_{n,k}$.  Then the above inequality $(\dagger\dagger)$ holds for all $n \ge 1$, $s \ge 0$ and all $k \ge 1$.  Let $v_1 = z$. We arrange these countably infinitely many points $c_{n,s,k}, \, n \ge 1$, $s \ge 0$, $k \ge 1$, in a sequence and call this sequence $v_2, v_3, v_4, \cdots$ (note that we have defined $v_1 = z$) and let $\mathcal S^{(1)}$, $\mathcal S^{(2)}$, $\mathcal S^{(3)}$, $\cdots$, be the resulting Cantor sets of totally transitive points of $f$, $\mathcal M$ the resulting infinite set of positive integers, $\mathbb S = \bigcup_{j = 1}^\infty S^{(j)}$, the resulting $\beta$-scrambled set (with respect to $\mathcal M$) of $f$ obtained by arguing as above with respect to the chosen set $\{ v_1, v_2, v_3, \cdots \}$.  

Let $c$ and $d$ be any (may be identical) points in $\mathbb S$ and let $t \ge s \ge 0$ be any integers. 

Since $v_1 = z$ is a fixed point of $f$ and since the set $\{ c, d \}$ is $\mathcal M$-synchronously proximal to $v_1 = z$, we easily obtain  that
\begin{equation*}
\liminf_{\substack{m \to \infty \\ m \in \mathcal M}} \rho\big(f^m(f^t(c)), f^m(f^s(d))\big) = \rho(z, z) = 0. \tag{$\ddagger\ddagger$}  
\end{equation*}
\indent On the other hand, suppose $t > s \ge 0$.  Since the set $\{ c, d \}$ is $\mathcal M$-synchronously proximal to $v_i$ for each $i \ge 1$, $\{ c, d \}$ is $\mathcal M$-synchronously proximal to $c_{n,s,k}$ for each $n \ge 1$, each $s \ge 0$, and each $k \ge 1$.  Consequently, for each $k \ge 1$, there exists a strictly increasing sequence $< r_{\ell,k} >_{\ell \ge 1}$, depending on $s$ and $t$, of positive integers in $\mathcal M$ such that $\lim_{\ell \to \infty} f^{r_{\ell,k}}(c) = \lim_{\ell \to \infty} f^{r_{\ell,k}}(d) = c_{t-s,s,k}$.  In particular, for each $k \ge 1$, we have \\
[10pt]\hspace*{.8in}$\limsup_{\substack{m \to \infty \\ m \in \mathcal M}} \rho\big(f^m(f^t(c)), f^m(f^s(d))\big) \ge \lim_{\ell \to \infty} \rho\big(f^{r_{\ell,k}}(f^t(c)), f^{r_{\ell,k}}(f^s(d))\big)$ \\
[5pt]\hspace*{1.2in}$= \lim_{\ell \to \infty} \rho\big(f^t(f^{r_{\ell,k}}(c)), f^s(f^{r_{\ell,k}}(d))\big)  
= \rho\big(f^t(c_{t-s,s,k}), f^s(c_{t-s,s,k})\big)$ \\
[5pt]\hspace*{1.6in}$= \rho\big(f^{t-s}(f^s(c_{t-s,s,k})), f^s(c_{t-s,s,k})\big) > \delta_{t-s} \, - \, \frac 1k \ge \delta \, - \, \frac 1k\qquad\qquad$ \\
[10pt]which implies that 
$\limsup_{\substack{m \to \infty \\ m \in \mathcal M}} \rho\big(f^m(f^t(c)), f^m(f^s(d))\big) \ge \delta$ when $t > s \ge 0$.  Furthermore, for $t = s$, $c \ne d$ and $\{ c, d \} \subset \mathbb S = \bigcup_{j=1}^\infty {\mathcal S}^{(j)}$, it follows from ($\dagger$) that 
$\limsup_{\substack{m \to \infty \\ m \in \mathcal M}} \rho\big(f^m(f^s(c)), f^m(f^s(d))\big) = \beta \ge \delta$.  These, combined with the above $(\ddagger\ddagger)$, imply that the set $\widehat {\mathbb S} = \bigcup_{i=0}^\infty \, f^i(\mathbb S)$ is a dense {\it invariant} $\delta$-scrambled set (with respect to $\mathcal M$) of totally transitive points of $f$ in $\mathcal X$.  Since, for each $n \ge 1$, $n$ divides all sufficiently large positive integers in $\mathcal M$, we obtain that $\widehat {\mathbb S}$ is a dense {\it invariant} $\delta$-scrambled set of $f^n$ for all $n \ge 1$.  This proves Part (5).
\hfill\sq

We say that $f$ is (topologically) mixing if, for any nonempty open sets $U$ and $V$, there exists a positive integer $N$ such that $f^k(U) \cap V \ne \emptyset$ for all $k \ge N$.  We note that continuous mixing maps have similar results as those for continuous weakly mixing maps except that 

\begin{itemize}
\item[(1)]
for continuous mixing maps, the scrambled sets ({\it invariant or not}) can be obtained with respect to any {\it given countably infinitely many} infinite sets of positive integers while those for continuous weakly mixing maps are with respect to {\it some} infinite sets of positive integers and 

\item[(2)]
the points $x_i$'s for continuous mixing maps $f$ have nonempty $\omega$-limit sets instead of having compact orbit closures $\overline{O_f(x_i)}$'s for continuous weakly mixing maps.  
\end{itemize}

Now suppose $f$ is mixing on $\mathcal X$ and $x_0$ is a point of $\mathcal X$ with an $\omega$-limit point $y_0$.  It is clear that
\begin{quote}
for any $n \ge 2$, any pairwise disjoint compact sets $K_1, K_2, \cdots, K_n$ with nonempty interiors, any nonempty open sets $V_1, V_2, \cdots, V_n$ in $\mathcal X$ and any open neighborhood $V(y_0)$ of $y_0$, there exists a positive integer $N$ such that $f^k\big(K_i\big) \cap V_i \ne \emptyset$ for all $1 \le i \le n$ and {\it all} $k > N$, and there exist countably infinitely many positive integers $k_j > N$ such that $f^{k_j}(x_0) \in V(y_0)$.  For each such $k_j > N$, we have $f^{k_j}\big(K_i\big) \cap V_i \ne \emptyset$ for all $1 \le i \le n$ and $f^{k_j}(x_0) \in V(y_0)$. \hfill (*)
\end{quote}

To prove the aforementioned similar results for continuous mixing maps, we can modify the proof of Theorem 4 with Lemmas 3(1) $\&$ 3(2) replaced by the above fact (*).  For convenience, for any infinite set $\mathcal M$ of positive integers and any positive integer $j$, we let $j\mathcal M = \{ jm: m \in \mathcal M \}$ and let $\mathcal M_1, \mathcal M_2, \mathcal M_3, \cdots$, be any sequence of infinite sets of positive integers.  Furthermore, we call the steps at the $\ell^{th}$ stage, for each $\ell \ge 1$, in the proof of Theorem 4 {\it original steps} and leave {\it the original first step} alone and call\\
[8pt]the group of the original second step to the original $(\ell+1)^{st}$ step {\it the original group 1}, \\
the group of the original $(\ell+2)^{nd}$ step to the original $(2\ell^2+3\ell+1)^{st}$ step {\it the original group 2}, \\
the group of all the original steps right after the original $(2\ell^2+3\ell+1)^{st}$ step {\it the original group 3}. \\
[8pt]\indent For each $\ell \ge 1$, at the $\ell^{th}$ stage, we expand each step of the $\ell$ steps in {\it the original group 1} into $\ell \cdot \ell$ steps by choosing, at each step, the respective $k_{\ell, \cdots}$ to be an (ever increasing) integer in $\mathcal M_1$ instead of being divisible by $\ell!$.  Then right after this modified group 1, we insert $\ell-1$ groups of steps each of which is similar to the just modified group 1 in the way that, for $2 \le j \le \ell$, in the $(j-1)^{st}$ inserted group, we choose, at each step, the respective $k_{\ell, \cdots}$ to be an (ever increasing) integer in $j\mathcal M_1$ instead of being divisible by $\ell!$.  So far, we obtain $\ell$ groups of steps, each group consists of $\ell$ steps.  Now, right after these $\ell$ groups of steps, we insert another similar $\ell$ groups of steps in the way that, for each $1 \le j \le \ell$, the respective integers $k_{\ell,\cdots}$'s in the $j^{th}$ group are all in $j\mathcal M_2$, and then followed by another $\ell$ groups of steps in the way that, for each $1 \le j \le \ell$, the respective integers $k_{\ell,\cdots}$'s in the $j^{th}$ group are all in $j\mathcal M_3$, $\cdots$, until the final $\ell$ groups of steps in the way that, for each $1 \le j \le \ell$, the respective integers $k_{\ell,\cdots}$'s in the $j^{th}$ group are all in $j\mathcal M_\ell$.    

So we have $\ell$ new groups of steps (each new group consists of $\ell \cdot \ell$ steps) followed by {\it the original group 2} which are related to the points $x_m$'s with nonempty $\omega$-limit sets.  If $x_0$ is a point in $\mathcal X$ with an $\omega$-limit point $y_0$, then, since every open neighborhood $V_0$ of $y_0$ contains $f^k(x_0)$ for infinitely many $k$'s, we modify {\it the original group 2} which are pertaining to the point $x_0$ by replacing Lemma 3(2) with the above fact (*) by letting $V_i = V_0$ for all $1 \le i \le n$ or $V_i = \widehat V_0$ for all $1 \le i \le n$, and replacing the open set $G_s(x_0)$ in Lemma 3(2) with $V_0$.  These new modified group 2 consists of $2\ell(\ell+1)$ steps.  See {\bf\cite{du4}} for details.

Finally, for each $\ell \ge 1$, at the $\ell^{th}$ stage, we modify {\it the original group 3} as we just did for {\it the original group 1} by choosing, at each step in the {\it the original group 3}, the respective $k_{\ell, \cdots}$ to be an (ever increasing) integer in $\mathcal M_1$ instead of being divisible by $\ell!$.  Then we insert successively $\ell-1$ groups of steps similar to the just modified group 3 such that, for $2 \le j \le \ell$, in the $(j-1)^{st}$ inserted group, we choose, at each step, the respective (ever increasing) integers $k_{\ell,\cdots}$ from the given set $j\mathcal M_1$.  So far, right after the modified group 2, we obtain \\
\centerline{$\ell$ {\it new} groups of steps each group consists of $(\ell+1)^{\ell \cdot 2^{\ell+1}}$ steps.}

We mimick these $\ell$ {\it new} groups of steps pertaining to $\mathcal M_1, 2\mathcal M_1, 3\mathcal M_1, \cdots, \ell\mathcal M_1$ respectively by replacing $\mathcal M_1$, at each step of these $\ell$ {\it new} groups, with $\mathcal M_2, \mathcal M_3, \cdots, \mathcal M_\ell$ successively and obtain \\
$\ell$ {\it new} groups of steps with all relevant integers $k_{\ell,\cdots}$'s in $\mathcal M_2, 2\mathcal M_2, 3\mathcal M_2, \cdots, \ell\mathcal M_2$ respectively and \\
$\ell$ {\it new} groups of steps with all relevant integers $k_{\ell,\cdots}$'s in $\mathcal M_3, 2\mathcal M_3, 3\mathcal M_3, \cdots, \ell\mathcal M_3$ respectively  \\
\centerline{and so on and so forth until we obtain} \\
$\ell$ {\it new} groups of steps with all relevant integers $k_{\ell,\cdots}$'s in $\mathcal M_\ell, 2\mathcal M_\ell, 3\mathcal M_\ell, \cdots, \ell\mathcal M_\ell$ respectively. 

In summary, at the $\ell^{th}$ stage, we have totally $1 + \ell \cdot (\ell \cdot \ell) + 2\ell(\ell+1) + \big((\ell+1)^{\ell \cdot 2^{\ell+1}}\big) \cdot (\ell \cdot \ell)$ steps.  The rest of the proof is similar to that of Theorem 4 with Lemma 3 replaced by the above (*).  We omit the details.  

\noindent
{\bf Theorem 5.}
{\it Let $\mathcal M_1, \mathcal M_2, \mathcal M_3, \cdots$, be a sequence of infinite sets of positive integers and let $(\mathcal X, \rho)$ be an infinite locally compact separable metric space with metric $\rho$ and let $f : \mathcal X \rightarrow \mathcal X$ be a continuous mixing map.  Let $\beta$ and $\delta$ be defined by putting
$$
\beta =  \begin{cases}
               \infty, & \text{if $\mathcal X$ is unbounded}, \cr
               diam(\mathcal X) = \sup\big\{ \rho(x, y): \{ x, y \} \subset \mathcal X \big\}, & \text{if $\mathcal X$ is bounded}, \cr
       \end{cases}
$$
and 
$$
\delta = \inf_{n \ge 1} \big\{ \sup \{ \rho\big(f^n(x), x\big): x \in \mathcal X \} \big\}.
$$ 

Let $\{ v_1, \, v_2, \, v_3, \, \cdots \}$ be any countably infinite subsets of $\mathcal X$ and let $\{ x_1, x_2, x_3, \cdots \}$ be a countably infinite set of points in $\mathcal X$ whose $\omega$-limit sets are nonempty.  Then $\delta > 0$ and there exist countably infinitely many pairwise disjoint Cantor sets ${\mathcal S}^{(1)}, \, {\mathcal S}^{(2)}, \, {\mathcal S}^{(3)}, \, \cdots$ of totally transitive points of $f$ in $\mathcal X$ such that 
\begin{itemize}
\item[{\rm (1)}]
for any integers $k \ge 1$, $\ell \ge 1$, $n \ge 1$, and any distinct points $a_1, a_2, \cdots, a_k$ in $\mathbb S = \bigcup_{j=1}^\infty {\mathcal S}^{(j)}$, the set $\{F_k^{nm_\ell}\big((a_1, a_2, \cdots, a_k)\big): m_\ell \in \mathcal M_\ell \}$ is dense in $\mathcal X_k = \mathcal X \times \mathcal X \times \cdots \times \mathcal X$ ($k$ terms), where $F_k: \mathcal X_k \longrightarrow \mathcal X_k$ is the map defined by $F_k\big((x_1, x_2, \cdots, x_k)\big)=\big(f(x_1), f(x_2), \cdots, f(x_k)\big)$.  Consequently, the point $(a_1, a_2, \cdots, a_k)$ is a totally transitive point of $F_k$ and, for each $n \ge 1$, the set $\{\big(F_k^n\big)^{m_\ell}\big((a_1, a_2, \cdots, a_k)\big): m_\ell \in \mathcal M_\ell \}$ is dense in $\mathcal X_k$;

\item[{\rm (2)}]
the set $\mathbb S = \bigcup_{j=1}^\infty {\mathcal S}^{(j)}$ \!\! is a dense $\beta$-scrambled set (with respect to $\mathcal M_\ell$ for all $\ell \ge 1$) of $f^n$ for all $n \ge 1$;

\item[{\rm (3)}]
for any integer $r \ge 1$, the set $\bigcup_{j=1}^r {\mathcal S}^{(j)}$ is $(n\mathcal M_\ell)$-synchronously proximal to $v_m$ for all $\ell \ge 1$, all $n \ge 1$ and all $m \ge 1$; 

\item[{\rm (4)}]
for any integers $r \ge 1$ and $s \ge 0$, the set $f^s\big(\bigcup_{j=1}^r {\mathcal S}^{(j)}\big)$ is dynamically synchronously proximal to $x_m$ for all $m \ge 1$ and, for any point $x$ in the set $\{x_1, x_2, x_3, \cdots\}$ and any point $c$ in the set $\widehat {\mathbb S} = \bigcup_{i=0}^\infty \, f^i({\mathbb S})$, we have $\liminf_{n \to \infty} \rho\big(f^n(x), f^n(c)\big) = 0$ and $\limsup_{n \to \infty} \rho\big(f^n(x), f^n(c)\big) \ge \beta/2$, i.e., $\{ x, c \}$ is a ($\beta/2$)-scrambled set of $f$;

\item[{\rm (5)}]
if $f$ has a fixed point, then the above Cantor sets ${\mathcal S}^{(1)}, {\mathcal S}^{(2)}, {\mathcal S}^{(3)}, \cdots$ of totally transitive points of $f$ can be chosen (by choosing appropriate $v_1, v_2, v_3, \cdots$) so that Parts (1), (2), (3) $\&$ (4) hold and, the set $\widehat {\mathbb S} = \bigcup_{i=0}^\infty \, f^i(\mathbb S)$ is a dense invariant $\delta$-scrambled set (with respect to $\mathcal M_\ell$ for all $\ell \ge 1$) of $f^n$ for all $n \ge 1$.
\end{itemize}}

Let $f$ be a continuous transitive map on an interval $I$.  It follows from {\bf\cite{mai}} that, for some $\va > 0$, $f$ has a dense $\va$-scrambled set which is a union of countably infinitely many synchronously proximal Cantor sets of transitive points of $f$.  For such maps, we have the following generalization.  

\noindent
{\bf Corollary 6.}
{\it Let $\mathcal M_1, \mathcal M_2, \mathcal M_3, \cdots$ be a sequence of infinite sets of positive integers.  Let $I$ be an interval in the real line and let $f : I \longrightarrow I$ be a continuous transitive map.  Then the following hold:

\begin{itemize}
\item[(1)]
if $f^2$ is transitive on $I$, then there exists a set $\mathbb S$ which is a union of countably infinitely many pairwise disjoint synchronously proximal Cantor sets of totally transitive points of $f$ such that $\mathbb S$ is a dense $(\sup I - \inf I)$-scrambled set (with respect to $\mathcal M_\ell$ for all $\ell \ge 1$) of $f^n$ in $I$ for all $n \ge 1$ and $\bigcup_{i=0}^\infty f^i(\mathbb S)$ is a dense invariant $\delta$-scrambled set (with respect to $\mathcal M_\ell$ for all $\ell \ge 1$) of $f^n$ in $I$ for all $n \ge 1$, where $\delta = \inf_{n \ge 1} \{ \sup_{x \in I} |f^n(x)-x|\} > 0$;

\item[(2)]
if $f^2$ is not transitive on $I$ and $\breve z$ is the unique fixed point of $f$ in $I$, then there exist a dense set $\mathbb S_1$ of totally transitive points of $f^2$ in $I \cap (-\infty, \breve z]$ and a dense set $\mathbb S_2$ of totally transitive points of $f^2$ in $I \cap [\breve z, \infty)$, both of which are unions of countably infinitely many pairwise disjoint synchronously proximal (with respect to $f^2$) Cantor sets, such that $\mathbb S_1 \cup \mathbb S_2$ is a dense $\be_{\min}$-scrambled set (with respect to $\mathcal M_\ell$ for all $\ell \ge 1$) of $f^{n}$ in $I$ for all $n \ge 1$ and, for each $k = 1, 2$, $\widehat {\mathbb S}_k \cup f(\widehat {\mathbb S}_k)$, where $\widehat {\mathbb S}_k = \bigcup_{i = 0}^\infty f^{2i}(\mathbb S_k)$, is a dense invariant $\delta_{\min}$-scrambled set (with respect to $\mathcal M_\ell$ for all $\ell \ge 1$) of $f^{n}$ in $I$ for all $n \ge 1$, where $\be_{\min} = \min \{ \breve z - \inf I, \, \sup I - \breve z \}$ and $\delta_{\min} = \min \big\{\inf_{n \ge 1} \{\sup_{x \in I \cap (-\infty, \breve z]} |f^{2n}(x)-x|\}$, $\inf_{n \ge 1} \{\sup_{x \in I \cap [\breve z, \infty)} |f^{2n}(x)-x|\}\big\} > 0$.
\end{itemize}}
\noindent
{\it Proof.}
If $f^2$ is transitive on $I$, then $f$ is mixing.  So, Part (1) follows from Theorem 5.  In the following, we assume that $f$ is transitive but not $f^2$ on $I$.  Then it is well-known that there exists a unique fixed point $\breve z$ of $f$ which lies in the interior of $I$ such that $f\big(I \cap (-\infty, \breve z]\big)$ is a dense subset of $I \cap [\breve z, \infty)$ and $f\big(I \cap [\breve z, \infty)\big)$ is a dense subset of $I \cap (-\infty, \breve z]$ and, on each of the intervals $I \cap (-\infty, \breve z]$ and $I \cap [\breve z, \infty)$, $f^2$ is transitive and has at least two fixed points (the unique fixed point $\breve z$ of $f$ and a period-2 point of $f$) and so is mixing.  It follows from Theorem 5 that both $\delta_1 = \inf_{n \ge 1} \big\{ \sup_{x \in I \cap (-\infty, \breve z]} |f^{2n}(x)-x| \big\}$ and $\delta_2 = \inf_{n \ge 1} \big\{ \sup_{x \in I \cap [\breve z, \infty)} |f^{2n}(x)-x| \big\}$ are positive and the following hold:  
\begin{quote}
for any $n \ge 2$, any pairwise disjoint nonempty compact intervals $K_1^{(1)}, K_2^{(1)}, \cdots, K_n^{(1)}$ 
and any nonempty open intervals $V_1^{(1)}, V_2^{(1)}, \cdots, V_n^{(1)}$ in $I \cap (-\infty, \breve z]$, and any pairwise 
disjoint nonempty compact intervals $K_1^{(2)}, K_2^{(2)}, \cdots, K_n^{(2)}$ and any nonempty open 
intervals $V_1^{(2)}, V_2^{(2)}, \cdots, V_n^{(2)}$ in $I \cap [\breve z, \infty)$, there exists a positive integer $N$ such that 
$f^k\big(K_i^{(j)}\big) \cap V_i^{(j)} \ne \emptyset$ for all $1 \le i \le n$, all $j = 1, 2$ and {\it all} $k \ge N$. \hfill{(**)}
\end{quote}
\indent By choosing $v_1^{(1)} = \breve z$ and appropriate $v_2^{(1)}, v_3^{(1)}, v_4^{(1)}, \cdots$ in $I \cap (-\infty. \breve z]$, and $v_1^{(2)} = \breve z$ and appropriate $v_2^{(2)}, v_3^{(2)}, v_4^{(2)}, \cdots$ in $I \cap [\breve z, \infty)$ and by mimicking the arguments in the proofs of Theorems 4 $\&$ 5 with Lemma 3(1) replaced by the above fact (**) to the map $f^2$ on both $I \cap (-\infty, \breve z]$ and $I \cap [\breve z, \infty)$ {\it synchronously}, we obtain that there exist dense sets $\mathbb S_1$ and $\mathbb S_2$ of {\it totally transitive points} of $f^2$ in $I \cap (-\infty, \breve z]$ and in $I \cap [\breve z, \infty)$ respectively both of which are unions of countably infinitely many pairwise disjoint synchronously proximal (with respect to $f^2$) Cantor sets such that 

\begin{itemize}
\item[{\rm (a)}]
for all $n \ge 1$, ${\mathbb S}_1$ is a dense $(\breve z - \inf I)$-scrambled set (with respect to $\mathcal M_\ell$ for all $\ell \ge 1$) of $(f^2)^n$ and $\widehat {\mathbb S}_1 = \bigcup_{i \ge 0} f^{2i}(\mathbb S_1)$ is a dense {\it invariant} $\delta_1$-scrambled set (with respect to $\mathcal M_\ell$ for all $\ell \ge 1$) of $(f^2)^n$ in $I \cap (-\infty, \breve z]$, ${\mathbb S}_2$ is a dense $(\sup I - \breve z)$-scrambled set (with respect to $\mathcal M_\ell$ for all $\ell \ge 1$) of $(f^2)^n$ and $\widehat {\mathbb S}_2 = \bigcup_{i \ge 0} f^{2i}(\mathbb S_2)$ is a dense {\it invariant} $\delta_2$-scrambled set (with respect to $\mathcal M_\ell$ for all $\ell \ge 1$) of $(f^2)^n$ in $I \cap [\breve z, \infty)$;

\item[{\rm (b)}]
for any integers $n \ge 1$ and $\ell \ge 1$ and, any two distinct points $u_1$ and $u_2$ in ${\mathbb S}_2$, the set $\{ \big(f^{2n} \times f^{2n}\big)^{m_\ell}\big((u_1, u_2)\big): m_\ell \in \mathcal M_\ell \}$ is dense in $\big(I \cap [\breve z, \infty)\big) \times \big(I \cap [\breve z, \infty)\big)$ and the set $\{ (f^{2n})^{m_\ell}(u_2): m_\ell \in \mathcal M_\ell \}$ is dense in $I \cap [\breve z, \infty)$ and

\item[{\rm (c)}]
for any integers $n \ge 1$ and $\ell \ge 1$ and, any two distinct points $u$ and $u'$ in $\mathbb S_1 \cup \mathbb S_2$, the set $\{ u, u' \}$ is $(n\mathcal M_\ell)$-synchronously proximal (with respect to $f^2$) to the unique fixed point $\breve z$ of $f$, and so, for any two distinct points $w$ and $w'$ in $\widehat {\mathbb S}_1 \cup f(\widehat {\mathbb S}_1) \cup \widehat {\mathbb S}_2 \cup f(\widehat {\mathbb S}_2)$, the set $\{ w, w' \}$ is $(n\mathcal M_\ell)$-synchronously proximal (with respect to $f$) to the fixed point $\breve z$ of $f$.  
\end{itemize}

In the following, we shall show that, $\mathbb S_1 \cup \mathbb S_2$ is a dense $\be_{\min}$-scrambled set (with respect to $\mathcal M_\ell$ for all $\ell \ge 1$) of $f^{n}$ in $I$ for all $n \ge 1$ and both $\widehat {\mathbb S}_1 \cup f(\widehat {\mathbb S}_1)$ and $\widehat {\mathbb S}_2 \cup f(\widehat {\mathbb S}_2)$ are dense {\it invariant} $\delta_{\min}$-scrambled set (with respect to $\mathcal M_\ell$ for all $\ell \ge 1$) of $f^{n}$ in $I$ for all $n \ge 1$.

Let $u$ and $u'$ be any two distinct points in $\mathbb S_1 \cup \mathbb S_2$.  

If both $u$ and $u'$ lie in the set $\mathbb S_1$ or the set $\mathbb S_2$, then it is clear that $\{ u, u' \}$ is a $\be_{\min}$-scrambled set (with respect to $\mathcal M_\ell$ for all $\ell \ge 1$) of $f^{2n}$ for all $n \ge 1$.  

Now, suppose $u \in \mathbb S_1$ and $u' \in \mathbb S_2$.  If follows from (c) above that, for each $n \ge 1$ and all $\ell \ge 1$, $\{ u, u' \}$ is $(n\mathcal M_\ell)$-synchronously proximal (with respect to $f^2$) to the fixed point $\breve z$ of $f$.  So, we obtain that
$0 \le \liminf_{\substack{m_\ell \to \infty \\ m_\ell \in \mathcal M_\ell}} \big|(f^{2n})^{m_\ell}(u) - (f^{2n})^{m_\ell}(u')\big| = \breve z - \breve z = 0$.  On the other hand, it is clear that, for all $i \ge 0$, $f^i(u)$ and $f^i(u')$ lie on opposite sides of $\breve z$.  Since, by (b) above, the set $\{ (f^{2n})^{m_\ell}(u'): m_\ell \in \mathcal M_\ell \}$ is dense in $I \cap [\breve z, \infty)$ for all $n \ge 1$ and $\ell \ge 1$, and since the orbit of $u$ with respect to $f^2$ (and hence to $f^{2n}$) is contained entirely in $I \cap (-\infty, \breve z]$, we easily obtain that $\limsup_{\substack{m_\ell \to \infty \\ m_\ell \in \mathcal M_\ell}} \big|(f^{2n})^{m_\ell}(u) - (f^{2n})^{m_\ell}(u')\big| \ge \sup I - \breve z \ge \be_{\min}$.  Consequently, $\{ u, u' \}$ is a $\be_{\min}$-scrambled set (with respect to $\mathcal M_\ell$ for all $\ell \ge 1$) of $f^{2n}$ for all $n \ge 1$.  

This confirms that $\mathbb S_1 \cup \mathbb S_2$ is a dense $\be_{\min}$-scrambled set (with respect to $\mathcal M_\ell$ for all $\ell \ge 1$) of $f^{2n}$ in $I$ for all $n \ge 1$.  In particular, $\mathbb S_1 \cup \mathbb S_2$ is a dense $\be_{\min}$-scrambled set (with respect to $\mathcal M_\ell$ for all $\ell \ge 1$) of $f^n$ in $I$ for all $n \ge 1$.  

We now show that the set $\widehat {\mathbb S}_2 \cup f(\widehat {\mathbb S}_2)$ \big($\widehat {\mathbb S}_1 \cup f(\widehat {\mathbb S}_1)$ respectively\big) is a dense {\it invariant} $\delta_{\min}$-scrambled set (with respect to $\mathcal M_\ell$ for all $\ell \ge 1$) of $f^n$ in $I$ for all $n \ge 1$. 

Let $w$ and $w'$ be any two distinct points in $\widehat {\mathbb S}_2 \cup f(\widehat {\mathbb S}_2)$.

If both $w$ and $w'$ lie in the set $\widehat {\mathbb S}_2$, then it is clear that $\{ w, w' \}$ is a $\delta_2$-scrambled set (with respect to $\mathcal M_\ell$ for all $\ell \ge 1$) of $f^{2n}$ for all $n \ge 1$, where $\delta_2 = \inf_{n \ge 1} \{\sup_{x \in I \cap [\breve z, \infty)} |f^{2n}(x)-x|\}\big\} > 0$.  

If both $w$ and $w'$ lie in the set $f(\widehat {\mathbb S}_2)$, then there exist two points $\hat w$ and $\hat w'$ in ${\mathbb S}_2$ and two {\it odd} integers $i \ge j \ge 1$ such that $f^i(\hat w) = w$ and $f^j(\hat w') = w'$, where either $\hat w \ne \hat w'$ or, $\hat w = \hat w'$ and $i > j \ge 1$.  We have two cases to consider:

Case 1.  $\hat w \ne \hat w'$ and $i \ge j \ge 1$.  In this case, it is clear that $f^j\big(I \cap [\breve z, \infty)\big)$ is dense in $I \cap (-\infty, \breve z]$.  Let $a_1, a_2, a_3, \cdots$ be a sequence of points in $I \cap (-\infty, \breve z)$ such that $\lim_{k \to \infty} a_k = \inf I$ and let $ <b_k >_{k \ge 1}$ be a sequence of points in $I \cap (\breve z, \infty)$ such that $f^j(b_k) = a_k$ for all $k \ge 1$.  Since, by (b) above, for each $n \ge 1$ and each $\ell \ge 1$, the set $\{(f^{2n} \times f^{2n})^{m_\ell}\big((\hat w, \hat w')\big): m_\ell \in \mathcal M_\ell \}$ is dense in $\hat I \times \hat I$, where $\hat I = I \cap [\breve z, \infty)$, we see that, for each $k \ge 1$ and each $\ell \ge 1$, there exist two sequences $< m_{\ell,p}>_{p \ge 1}$ and $< m_{\ell,k,p} >_{p \ge 1}$ of integers in $\mathcal M_\ell$ such that $\lim_{p \to \infty} (f^{2n} \times f^{2n})^{m_{\ell,p}}\big((\hat w, \hat w')\big) = (\breve z, \breve z)$ and $\lim_{p \to \infty} (f^{2n} \times f^{2n})^{m_{\ell,k,p}}\big((\hat w, \hat w')\big) = (\breve z, b_k)$.  Consequently, since $\lim_{p \to \infty} (f^{2n} \times f^{2n})^{m_{\ell,p}}(w, w') = \lim_{p \to \infty} (f^{2n} \times f^{2n})^{m_{\ell,p}}\big(f^i(\hat w), f^j(\hat w')\big) = \big(f^i(\breve z), f^j(\breve z)\big) = (\breve z, \breve z)$ and since $\lim_{p \to \infty} (f^{2n} \times f^{2n})^{m_{\ell,k,p}}(w, w') = \lim_{p \to \infty} (f^{2n} \times f^{2n})^{m_{\ell,k,p}}\big(f^i(\hat w), f^j(\hat w')\big) = \big(f^i(\breve z), f^j(b_k)\big) = (\breve z, a_k)$ and $\lim_{k \to \infty} a_k = \inf I$, by taking the sequences $<m_{\ell,p}>_{p \ge 1}$ and $<m_{\ell,p,p}>_{p \ge 1}$, we obtain that $\{ w, w' \}$ is a $(\breve z - \inf I)$-scrambled set (with respect to $\mathcal M_\ell$ for all $\ell \ge 1$) of $f^{2n}$ for all $n \ge 1$.  

Case 2.  $\hat w = \hat w'$ and $i > j \ge 1$.  In this case, since both $i$ and $j$ are odd integers, $i-j$ is even.  For any $k \ge 1$, let $a_k$ be a point in $I \cap (-\infty, \breve z)$ such that $|f^{i-j}(a_k) - a_k| > \big(\sup_{x \in I \cap (-\infty, \breve z]} |f^{i-j}(x)-x|\big) - 1/k$ and let $b_k$ be a point in $I \cap (\breve z, \infty)$ such that $f^j(b_k) = a_k$.  Since, by (b) above, for each $n \ge 1$ and each $\ell \ge 1$, the set $\{(f^{2n})^{m_\ell}(\hat w): m_\ell \in \mathcal M_\ell \}$ is dense in $I \cap [\breve z, \infty)$, we see that, for each $k \ge 1$ and each $\ell \ge 1$, there exist two sequences $< m_{\ell,p}'>_{p \ge 1}$ and $< m_{\ell,k,p}'>_{p \ge 1}$ of integers in $\mathcal M_\ell$ such that $\lim_{p \to \infty} (f^{2n})^{m_{\ell,p}'}(\hat w) = \breve z$ and $\lim_{p \to \infty} (f^{2n})^{m_{\ell,k,p}'}(\hat w) = b_k$.  Consequently, we obtain that $0 \le \liminf_{\substack {m_\ell \to \infty \\ m_\ell \in \mathcal M_\ell}} |(f^{2n})^{m_\ell}(w)-(f^{2n})^{m_\ell}(w')| \le \lim_{p \to \infty} |(f^{2n})^{m_{\ell,p}'}\big(f^i(\hat w)\big)-(f^{2n})^{m_{\ell,p}'}\big(f^j(\hat w)\big)| = f^{i}(\breve z) - f^{j}(\breve z) = 0$ and $\limsup_{\substack {m_\ell \to \infty \\ m_\ell \in \mathcal M_\ell}} |(f^{2n})^{m_\ell}(w)-(f^{2n})^{m_\ell}(w')| \ge 
\lim_{p \to \infty} |(f^{2n})^{m_{\ell,k,p}'}(w)- (f^{2n})^{m_{\ell,k,p}'}(w')| = |f^{i}(b_k)-f^{j}(b_k)| = |f^{i-j}\big(f^j(b_k)\big)-f^{j}(b_k)| = |f^{i-j}(a_k)-a_k| > \big(\sup_{x \in I \cap (-\infty, \breve z]} |f^{i-j}(x)-x|\big) - 1/k$ for all $k \ge 1$.  By letting $k$ tend to $\infty$, we obtain that 
$\limsup_{m_\ell \to \infty} |(f^{2n})^{m_\ell}(w)-(f^{2n})^{m_\ell}(w')| \ge \sup_{x \in I \cap (-\infty, \breve z]} |f^{i-j}(x)-x| \ge \delta_1$, where $\delta_1 = \inf_{\hat n \ge 1} \{\sup_{x \in I \cap (-\infty, \breve z]} |f^{2\hat n}(x)-x|\} > 0$.  This implies that $\{ w, w' \}$ is a $\delta_1$-scrambled set (with respect to $\mathcal M_\ell$ for all $\ell \ge 1$) of $f^{2n}$ for all $n \ge 1$. 

Now suppose $w \in \widehat {\mathbb S}_2$ and $w' \in f(\widehat {\mathbb S}_2)$.  Let $\tilde w$ and $\tilde w'$ be points in ${\mathbb S}_2$ such that $f^k(\tilde w) = w$ and $f^j(\tilde w') = w'$ for some {\it even} $k \ge 0$ and some {\it odd} $j \ge 1$.  Since both $\tilde w$ and $\tilde w'$ lie in ${\mathbb S}_2$, it follows from (c) above that, for each $n \ge 1$ and each $\ell \ge 1$, $\big\{ \tilde w, \tilde w' \big\}$ is $(n\mathcal M_\ell)$-synchronously proximal (with respect to $f^2$) to $\breve z$.  Since $\breve z$ is a fixed point of $f$, we easily obtain that $\liminf_{\substack{m_\ell \to \infty \\ m_\ell \in \mathcal M_\ell}} \big|(f^{2n})^{m_\ell}(w) - (f^{2n})^{m_\ell}(w')\big| = 0$.  On the other hand, it is clear that, for each $i \ge 0$, $f^i(w)$ and $f^i(w')$ lie on opposite sides of $\breve z$.  Since, by (b) above, for each $n \ge 1$ and each $\ell \ge 1$, the set $\{ (f^{2n})^{m_\ell}(w): m_\ell \in \mathcal M_\ell \}$ is dense in $I \cap [\breve z, \infty)$ and since the orbit $O_{f^{2n}}(w')$ of $w'$ with respect to $f^{2n}$ is contained entirely in $I \cap (-\infty, \breve z]$, we easily obtain that 
$\limsup_{\substack{m_\ell \to \infty \\ m_\ell \in \mathcal M_\ell}}\big|(f^2)^{nm_\ell}(w) - (f^2)^{nm_\ell}(w')\big| \ge \breve z - \inf I \ge \delta_{\min}$.  This implies that $\{ w, w' \}$ is a $\delta_{\min}$-scrambled set (with respect to $\mathcal M_\ell$ for all $\ell \ge 1$) of $f^{2n}$ for all $n \ge 1$.

Therefore, $\widehat {\mathbb S}_2 \cup f(\widehat {\mathbb S}_2)$ \big(and similarly $\widehat {\mathbb S}_1 \cup f(\widehat {\mathbb S}_1)$\big) is a dense {\it invariant} $\delta_{\min}$-scrambled set (with respect to $\mathcal M_\ell$ for all $\ell \ge 1$) of $f^{2n}$ in $I$ for all $n \ge 1$.  Since $f(\widehat {\mathbb S}_2 \cup f(\widehat {\mathbb S}_2))\! \subset\! f(\widehat {\mathbb S}_2) \cup f^2(\widehat {\mathbb S}_2) \subset f(\widehat {\mathbb S}_2) \cup \widehat {\mathbb S}_2$, we obtain that $\widehat {\mathbb S}_2 \cup f(\widehat {\mathbb S}_2)$ \big(and similarly $\widehat {\mathbb S}_1 \cup f(\widehat {\mathbb S}_1)$\big) is also a dense {\it invariant} $\delta_{\min}$-scrambled set (with respect to $\mathcal M_\ell$ for all $\ell \ge 1$) of $f^n$ in $I$ for all $n \ge 1$.  
\hfill\sq

\noindent
{\bf Remark.}
Let $(\mathcal X, \rho)$ be an infinite metric space with metric $\rho$ and let $f : \mathcal X \longrightarrow \mathcal X$ be a continuous map. 
In {\bf\cite{du3,du7,du4,du5,du6,v2,v1}} (see also {\bf\cite{wang}}), we say that $f$ is a chaotic map if there exists a positive number $\eta$ such that, for any point $x$ and any nonempty open set $V$ (not necessarily an open neighborhood of $x$) in $\mathcal X$, there exists a point $y$ in $V$ such that 
$$
\liminf_{n \to \infty} \rho\big(f^n(x), f^n(y)\big) = 0 \,\,\, \text{and} \,\,\, \limsup_{n \to \infty} \rho\big(f^n(x), f^n(y)\big) \ge \eta.
$$
In light of Theorem 4, we see that if $\mathcal X$ is an infinite compact metric space (without isolated points), then every continuous weakly mixing map $f$ from $\mathcal X$ into itself is {\it chaotic}.  However, not all {\it chaotic} maps are weakly mixing or transitive.  The following map $g : [0, 3] \longrightarrow [0, 3]$ is such an example: $g(x) = 2x$ for $0 \le x \le 1$; $g(x) = 4-2x$ for $1 \le x \le 2$ and $g(x) = 2x-4$ for $2 \le x \le 3$.  We refer to {\bf\cite{du5}} for details.  

\noindent
{\bf Acknowledgment}\\
This work was partially supported by the Ministry of Science and Technology of Taiwan.

\bibliographystyle{plain}

\end{document}